\documentclass[a4paper,reqno]{amsart}

\title{Operads and Jet modules}
\author{Marc A.~Nieper-Wi\ss kirchen}
\date{\today}
\address{Institut f\"ur Mathematik \\ Johannes Gutenberg-Universit\"at \\
  Staudinger Weg 9 \\ 55128 Mainz \\ Germany}
\email{nieper@mathemtik.uni-mainz.de}

\usepackage{amscd}
\usepackage{amssymb}
\usepackage{stmaryrd}
\usepackage[all,cmtip,poly]{xy}
\usepackage{graphicx}

\theoremstyle{plain}
\newtheorem{theorem}{Theorem}
\newtheorem{proposition}{Proposition}
\newtheorem{lemma}{Lemma}

\theoremstyle{definition}
\newtheorem{definition}{Definition}
\theoremstyle{remark}
\newtheorem{remark}{Remark}
\newtheorem{example}{Example}

\newcommand{\card}[1]{\left|#1\right|}

\newcommand{\cofin}{\operatorname{cofin}}
\newcommand{\colim}{\varinjlim}

\newcommand{\diff}{\bar\partial}
\newcommand{\Ext}{\mathrm{Ext}}
\newcommand{\id}{\mathrm{id}}

\newcommand{\SG}{\mathfrak S}
\newcommand{\cat}[1]{\mathcal #1}
\newcommand{\set}[1]{\mathbf #1}
\newcommand{\unit}{\mathbf 1}

\newcommand{\sgm}[1]{\mathcal #1}
\newcommand{\Ass}{\sgm{Ass}}
\newcommand{\Com}{\sgm{Com}}
\newcommand{\End}{\sgm{End}}
\newcommand{\ad}{\operatorname{ad}}
\newcommand{\op}{\circ}
\newcommand{\Ho}{\operatorname{Ho}}
\newcommand{\cone}{\operatorname{cone}}
\newcommand{\freeto}{\to_0}
\newcommand{\stableto}{\rightharpoonup}
\newcommand{\freestableto}{\rightharpoonup_0}
\newcommand{\im}{\operatorname{im}}

\begin{document}

\begin{abstract}
  Let $A$ be an algebra over an operad in a cocomplete closed
  symmetric monoidal category. We study the category of $A$-modules.
  We define certain symmetric product functors of such modules
  generalising the tensor product of modules over commutative
  algebras, which we use to define the notion of a jet module. This in
  turn generalises the notion of a jet module over a module over a
  classical commutative algebra.  We are able to define Atiyah classes
  (i.e.~obstructions to the existence of connections) in this
  generalised context. We use certain model structures on the category
  of $A$-modules to study the properties of these Atiyah classes.
  
  The purpose of the paper is not to present any really deep theorem.
  It is more about the right concepts when dealing with modules over
  an algebra that is defined over an arbitrary operad, i.e.~the aim is
  to show how to generalise various classical constructions, including
  modules of jets, the Atiyah class and the curvature, to the operadic
  context.
  
  For convenience of the reader and for the purpose of defining the
  notations, the basic definitions of the theory of operads and model
  categories are included.
\end{abstract}

\maketitle

\tableofcontents

\section*{Introduction}

Let $X$ be a complex manifold. We denote its $\mathcal O_X$-module of
K\"ahler differentials by $\Omega_X$. Let $\mathcal E$ be another
locally free $\mathcal O_X$-module. A (global) holomorphic connection
on $\mathcal E$ is a $\set C$-linear morphism $\nabla: \mathcal E \to
\Omega_X \otimes \mathcal E$ of sheaves on $X$ such that
\begin{equation}
  \label{eq:leibniz_for_connections}
  \nabla(f s) = df \otimes s + f \nabla s
\end{equation}
for a local section $f$ of $\mathcal O_X$ and a local section $s$ of
$\mathcal E$.  In general, no global holomorphic connection does
exist. The obstruction to the existence is given by the so-called
\emph{Atiyah class $\alpha_{\mathcal E} \in \Ext_X^1(\mathcal E,
  \Omega_X \otimes \mathcal E)$ of $\mathcal E$} (\cite{Atiyah57}).
The Atiyah class gives us thus cohomological invariants of the module
$\mathcal E$. In fact, the Chern classes of $\mathcal E$ can be
calculated from its Atiyah class. For details, we refer the reader to
the first secion of \cite{Kapranov99} or to the first chapter of our
book~\cite{Nieper04}.

Let us denote by $\mathcal A_X^{p, q}$ the $\mathcal C^\infty_X$-sheaf
of the $(p, q)$-differential forms on $X$. A $\mathcal
C^\infty$-connection on $\mathcal E$ is a $\set C$-linear morphism
$\nabla: \mathcal E \to \mathcal A^{1, 0}_X \otimes \mathcal E$ of
sheaves on $X$ such that~\eqref{eq:leibniz_for_connections} holds.
Using a partition of unity, one can show that such a connection always
exists for $\mathcal E$. As shown in, e.g.,~\cite{Kapranov99}, the
Atiyah class of $\mathcal E$ is then given by extension class
represented by $- \diff \nabla$, which is in fact $\mathcal
O_X$-linear (the occurrance of the sign depends on the identification
of the Dolbeault cohomology groups $H_{\bar \partial}^*(X, \Omega_X
\otimes \mathcal E^\vee \otimes \mathcal E)$ of $\Omega_X \otimes
\mathcal E^\vee \otimes \mathcal E$ with the Ext-groups
$\Ext^1_X(\mathcal E, \Omega_X \otimes \mathcal E)$ of $\mathcal E$ by
$\Omega_X \otimes \mathcal E$), see~\cite{Demailly97}.

In this article, we want to generalise these constructions
considerably. First of all, the notion of a connection does not truely depend
on the universal derivation $d: \mathcal O_X \to \Omega_X$. In fact, one may
substitute $d: \mathcal O_X \to \Omega_X$ by any derivation $d: \mathcal O_X
\to \mathcal M$ where $\mathcal M$ is a locally free $\mathcal
O_X$-module.

In~\cite{Kapranov99} another way besides the Dolbeault method to
calculate the Atiyah is presented: one may use \v Cech resolutions
instead of Dolbeault resolution, while substituting the Dolbeault
differential with the \v Cech differential. Formally, both methods are
very similar. The idea to handle not only these two resolutions
uniformly is to think of these resolutions as fibrant replacements
each with respect to a certain model structures on the category of
cochain complexes of $\mathcal O_X$-modules, where the weak
equivalences are the quasi-isomorphisms. The assumption that $\mathcal
M$ and $\mathcal E$ are locally free translate into the assumption
that both are cofibrant modules. Thus the theory of model categories
enters at this point and generalises, e.g.,~the \v Cech construction
mentioned above.

The next generalisation is to replace the category of $\mathcal O_X$-modules
by any closed symmetric monoidal abelian category whose category of cochain
complexes $\cat C^*$ possesses a suitable model structure. Given a commutative
algebra $A$ in this category $\cat C^*$, we can talk about $A$-modules. Thus
we are lead to the notion of an Atiyah class for (cofibrant) $A$-modules and
can try to translate everything else from the theory on complex
manifolds to this context.

The last step in generalising these constructions and notions is the
main reason why this article has been written. Instead of saying that
$A$ is a commutative algebra in $\cat C^*$, we may say that $A$ is an
algebra over the operad $\Com$ in $\cat C^*$
(see~\cite{MarklShniderStasheff02}). Replacing the operad $\Com$ by
any other operad, we arrive at other types of algebras, say
associative or Lie algebras. One can again talk of $A$-modules and we
show what to do to get the notion of connections and the Atiyah class
also in this operadic context.

There are still ways to generalise the results in this paper even
more.  Firstly, we only consider the case of cofibrant objects when
dealing with their Atiyah class. Similar, the module of differentials
is assumed to be cofibrant. One may try to extend the theory to the
non-cofibrant case.  Secondly, the theory developed here works for
model categories of cochain complexes, i.e.~for model categories being
modules over the category of cochain complexes of abelian groups. To
deal with the non-abelian case, one may work with model categories
that are modules over the category of simplicial sets. To deal with
this case one has to extend the notion of a derivation and a
connection to the simplicial case, and one can hope to get the notion
of an Atiyah class as well in this case. We haven't pursued this way
in this article.

\subsection{Organisation}
The article is organised in sections (numbered by 1., 2., etc.), which
themselves are divided into subsections (numbered by 1.1., 1.2., etc.). The
outline of the sections is as follows:

The first section deals with the notion of a symmetric closed monoidal
category, i.e.~with a categorical generalisation of the presence of the tensor
product in the category of modules over a commutative ring with a unit. We do
this so that we can talk about monoids (associative algebras) in these
categories and (left) modules over them. We shall apply these notions to
categories of cochain complexes.

The next section, which is on operadic notions, brings together all
the results of the theory of operads needed here. The first
subsections contain material well-known to people working on operads.
The other subsections deal with the generalisation of the theory of
connections of modules over a commutative algebra $A$ to modules over
arbitrary algebras. A main point here is to generalise the tensor
product over $A$ between modules. The answer we have come with is
given by what we call the \emph{lax product}. It gives back the tensor
product in the commutative case, which gives the category of
$A$-algebras a structure of a symmetric monoidal category. In general
however, our lax product only gives a structure of a so-called lax
symmetric monoidal category (\cite{DayStreet03}). Thus the name given
to our product.

Having the application to the category of cochain complexes of $\mathcal
O_X$-modules over a complex manifold $X$ in mind, we develop the necessary
portions of the model category theory to deal with suitable model structures
on the categories of modules over arbitrary algebras in the third and forth
section. These sections are mainly for the convenience of the reader and to
fix the needed notions. For any further study of model categories, we highly
recommend the book~\cite{Hovey99} and the references therein.

In the last section, we finally deal with the Atiyah class in the
operadic case. The Atiyah class becomes a morphism in a homotopy
category (which is nothing else than a derived category as our weak
equivalences are always exactly the quasiisomorphisms) and we show how
it can be calculated by using fibrant resolutions.
In~\cite{Kapranov99}, cohomological Bianchi identities for the Atiyah
class are proven. We prove these identities in the more general
operadic case. We also talk about the notion of (higher) curvature
classes. Again we have been inspired by the methods and results
in~\cite{Kapranov99}. Our final application is as follows: Starting
with a deformation of a free algebra $A$, we consider its category of
modules and show that the Atiyah class of its module $M$ of K\"ahler
differentials on $A$ defines a deformation of the free $A$-algebra
over $M$.  This gives a map compatible with gauge equivalence from the
solution set of one Maurer--Cartan equation to the solution set of
another Maurer--Cartan equation. In particular, every Lie algebra
$\mathfrak g$ (in fact every stong homotopy Lie algebra) gives rise to
a deformed version of the free commutative algebra over the vector
space $\mathfrak g$, over which we can study our geometrical notions
like jet modules and curvature forms. (This has been inspired
by~\cite{Kontsevich99}.)

The article includes a short appendix in which we give a proof for the
result that every object in a Grothendieck category $\cat A$ is small.
This is needed to use Quillen's small object argument
(see~\cite{Hovey99}) in order to put a model structure on the category
of cochain complexes over $\cat A$, which we want to do. The result
itself is well-known and one of the different types of proofs in the
literature can be found in~\cite{Hovey01}. However, all published
proofs we know use deep theorems about Grothendieck categories (like
the embedding theorem). Therefore we feet that it is time to give a
simple proof using just the basic properties of Grothendieck
categories.

We hope that our developed notions will prove their usefulness in
applications of the theory of modules over arbitrary algebras.

\subsection{A few remarks}

A lot of material in this article is included to round up the whole
exposition but is otherwise well-known to people working in the
particular fields. In these cases, we often give references to the
literature. However, these are often not the original works where the
material originally comes from but textbooks which may be better
accessible to the average reader.

A lot of the propositions in this article just end with the ``proof end
symbol'', meaning that no proof is included. This usually means that the proof
is straight-forward and follows directly from the definitions. Most of the
morphisms used in this paper are defined by some diagrams in some categories.
In order to keep the diagrams simple, a lot of arrows between two objects, say
$X$ and $Y$, are not annotated. This means that the arrow stands for the
``most natural'' morphism between the two objects. For example, an arrow $R
\otimes_k M \to M$ where $R$ is a $k$-algebra, $k$ a field, and $M$ is an
$R$-module over $k$ stands for the operation of $R$ on $M$, i.e.~the scalar
multiplication.

\section{Symmetric closed monoidal categories}

For the convenience of the reader, we recall in this section the
notion of a symmetric closed monoidal category, on which the whole
theory of operads is based.

\subsection{Tensor products}

Let $\cat C$ be a category. By a \emph{tensor product on $\cat C$}, we
understand a bifunctor $\cdot \otimes \cdot: \cat C \times \cat C \to
\cat C$ that is suitably associative, i.e.~there are fixed
isomorphisms, called \emph{associators}, $X \otimes (Y \otimes Z) \to
(X \otimes Y) \otimes Z$ natural in $X$, $Y$ and $Z$ for which the
\emph{coherence diagrams}
\[\xymatrix{
  & X \otimes (Y \otimes (Z \otimes W)) \ar[dr] \ar[dl] \\
  X \otimes ((Y \otimes Z)) \otimes W) \ar[d] & &
  (X \otimes Y) \otimes (Z \otimes W) \ar[d]
  \\
  (X \otimes (Y \otimes Z)) \otimes W \ar[rr] & &
  ((X \otimes Y) \otimes Z) \otimes W
}\]
that are built up from the associators and are natural in $X$, $Y$,
$Z$ and $W$, commute. 

\begin{remark}
  By Mac Lane's coherence theorem (\cite{MacLane98}), the commutativity of
  the pentagons above suffices to show that in fact all natural diagrams built
  up from the associators do commute.
\end{remark}

\subsection{Units}

Let $\otimes$ be a tensor product on $\cat C$. By a \emph{unit for $\otimes$}
we understand an object $\unit$ of $\cat C$ with fixed isomorphisms, called
the \emph{left} and the \emph{right unit law}, $\unit \otimes X \to X$ and $X
\otimes \unit \to X$, natural in $X$, such that the \emph{coherence diagrams}
\[\xymatrix{
  X \otimes (\unit \otimes Y) \ar[rr] \ar[dr] & & (X \otimes
  \unit) \otimes Y \ar[dl] \\
  & X \otimes Y,
}\]
that are built up from the associators and the unit laws and are
natural in $X$ and $Y$, commute. 

\begin{remark}
  Again, the commutativity of the triangles above (together with the
  commutativity of the pentagons above) suffices to make all natural
  diagrams built up from the associators and unit laws commute. 
\end{remark}

\begin{definition}
  A category together with a tensor product and a unit as defined
  above is a \emph{monoidal category}. 
\end{definition}

\subsection{Symmetric tensor products}

Let $\cat C$ be a monoidal category. The tensor product $\otimes$ is
called \emph{symmetric} if there are fixed isomorphisms, called
\emph{(symmetric) braidings}, $\gamma: X \otimes Y \to Y \otimes X$
natural in $X$ and $Y$ for which the \emph{coherence diagrams}
\[\xymatrix{
  X \otimes Y \ar[rr] \ar@{=}[dr] & & Y \otimes X \ar[dl] \\
  & X \otimes Y,
}\]
and
\[\xymatrix{
  & X \otimes (Y \otimes Z) \ar[dl] \ar[dr] \\
  (X \otimes Y) \otimes Z \ar[d] & &
  X \otimes (Z \otimes Y) \ar[d] \\
  Z \otimes (X \otimes Y) \ar[dr] & &
  (X \otimes Z) \otimes Y  \\
  & (Z \otimes X) \otimes Y \ar[ur]
}\]
that are built from the associators and the braidings and are natural
in $X$, $Y$ and $Z$, commute. 

\begin{remark}
  The commutativity of all coherence diagrams suffices to make all
  natural diagrams built up from the associators, unit laws and
  braidings commute. 

  Thus, a specific bracketing or ordering in iterated tensor products
  does not matter up to a uniquely defined natural isomorphism. 
\end{remark}

\begin{definition}
  A monoidal category with a symmetric tensor product as defined above
  is a \emph{symmetric monoidal category}. 
\end{definition}

\subsection{Symmetric monoidal functors}

Later we shall use the notion of a \emph{symmetric monoidal functor},
which is more or less a functor between symmetric monoidal categories
that respects the symmetric monoidal structures. We
follow~\cite{BergerMoerdijk03}. The authors of this article in turn
refer to~\cite{MandellMaySchwedeShipley01}.

Let $\cat C$ and $\cat D$ be two symmetric monoidal categories. Let
$F: \cat C \to \cat D$ be a functor. Assume $F$ comes equipped with
a morphism $\unit \to F(\unit)$ in $\cat D$ and morphisms
$F(X) \otimes F(Y) \to F(X \otimes Y)$ in $\cat D$ that are natural in
$X$ and $Y$.

The functor $F$ is \emph{(left) unital} if the diagrams
\[\xymatrix{
  \unit \otimes F(X) \ar[r]\ar[d] & F(X) \ar@{<-}[d]
  \\
  F(X) \otimes F(X) \ar[r] & F(\unit \otimes X),
}\]
natural in $X$, commute.

The functor $F$ is \emph{symmetric} if the diagrams
\[\xymatrix{
  F(X) \otimes F(Y) \ar[r]\ar[d] & F(Y) \otimes F(X) \ar[d] \\
  F(X \otimes Y) \ar[r] & F(Y \otimes X),
}\]
natural in $X$ and $Y$, commute.

Finally, the functor $F$ is \emph{associative} if the diagrams
\[\xymatrix{
  & F(X) \otimes (F(Y) \otimes F(Z)) \ar[dr]\ar[dl] \\
  (F(X) \otimes F(Y)) \otimes F(Z) \ar[d] & &
  F(X) \otimes F(Y \otimes Z) \ar[d] \\
  F(X \otimes Y) \otimes F(Z) \ar[dr] & &
  F(X \otimes (Y \otimes Z)) \ar[dl] \\
  & F((X \otimes Y) \otimes Z)
}\]
natural in $X$, $Y$ and $Z$, commute. 

\begin{definition}
  A functor $F$ as above that comes equipped with the morphism $\unit
  \to F(\unit)$ and the morphisms $F(X) \otimes F(Y) \to F(X \otimes
  Y)$, natural in $X$ and $Y$, which is unital, symmetric, and
  associative is a \emph{symmetric monoidal functor}.
\end{definition}

\subsection{Inner hom's}

Let $\cat C$ be a symmetric monoidal category. An \emph{inner hom for
  $\cat C$} is a bifunctor $\hom: \cat C \times \cat C \to \cat C$
such that there are bijections $\cat C(X \otimes Y, Z) \to \cat C(X,
\hom(Y, Z))$ natural in $X$, $Y$ and $Z$. In other words, $\hom(Y,
\cdot): \cat C \to \cat C$ is right adjoint to $\cdot \otimes Y: \cat
C \to \cat C$ for each object $Y \in \cat C$. 

\begin{definition}
  A \emph{closed symmetric monoidal category} is a symmetric
  monoidal category that possesses an inner hom as defined above. 
\end{definition}

\begin{example}
  One of the most prominent example of such a category is the category
  of $k$-modules for a commutative ring with unit $k$. The tensor
  product is the tensor product of $k$-modules over $k$, the unit is
  the object $k$. The inner hom of two objects is the set of
  $k$-linear maps between these objects considered as a $k$-module. 
  Finally, the associators, unit laws, symmetric braidings and
  adjunctions are given by the obvious isomorphisms. 
\end{example}

\begin{remark}
  Let us remark that in a closed symmetric monoidal category the
  functor $\cdot \otimes Y: \cat C \to \cat C$ does commute with every
  colimit for each object $Y$. This relies solely on the fact that the
  functor possesses a right adjoint. 
\end{remark}

Let us end this paragraph with a definition concerning the previously
defined structures in the case that the underlying category is
additive.
\begin{definition}
  A \emph{closed symmetric monoidal additive category} is a symmetric
  monoidal category that is at the same time an additive category and
  such that the symmetric monoidal category structures are
  compatible with the group laws on the hom-sets. 
\end{definition}

\subsection{Cochain complexes}

Let $\cat C$ be an additive category. By a \emph{cochain complex over
  $\cat C$} we understand a family $(X^n)_{n \in \mathbf Z}$ of
objects $X^n \in \cat C$ together with morphisms $\diff^n: X^n \to
X^{n + 1}$ with $\diff^{n + 1} \circ \diff^n = 0$ for all $n \in
\mathbf Z$. Such a family is usually abbreviated by $X^*$ or simply by
$X$. The morphisms $\diff^n$ are the \emph{differentials of $X^*$}. A
\emph{morphism between cochain complexes $X^*$ and $Y^*$ over $\cat
  C$} is a family $(f^n: X^n \to Y^n)_{n \in \mathbf Z}$ of morphisms
in $\cat C$ such that $\diff^n \circ f^n - f^{n + 1} \circ \diff^n =
0$ for all $n \in \mathbf Z$. One can \emph{concatenate} morphisms of
cochain complexes over $\cat C$ componentwise in the obvious way. This
leads to the following definition.
\begin{definition}
  The \emph{category of cochain complexes over $\cat C$} is the
  category whose objects are the cochain complexes over $\cat
  C$ and whose morphisms are the morphisms between cochain complexes
  defined above. It is denoted by $\cat C^*$. 
\end{definition}

For each cochain complex $X^*$ and each $n \in \set Z$ we define the
$n$-th cohomology group $H^n(X^*) := \ker \diff^n/\im \diff^{n -
  1}$. Note that each morphism $f: X^* \to Y^*$ of cochain complexes
induces as usual group homomorphisms between the cohomology groups of
$X^*$ and $Y^*$. We need this for the following (well-known)
definition:
\begin{definition}
  The morphism $f: X^* \to Y^*$ is a \emph{quasiisomorphism} if $f^*:
  H^*(X^*) \to H^*(Y^*)$ is an isomorphism on each cohomology group.
\end{definition}

We shall also make use of the notion of a \emph{free morphism between
  cochain complexes over $\cat C$}. A free morphism between two
cochain complexes $X^*$ and $Y^*$ is just a family $(f^n: X^n \to
Y^n)_{n \in \set Z}$ of morphisms in $\cat C$ with no compatibility
condition with the differentials whatsoever. Being a free morphism is
denoted by $f: X \freeto Y$.

Assume that $\cat C$ is a bicomplete (i.e.~arbitrary small limits and
colimits do exist) closed symmetric monoidal additive category. In
what follows we describe how to put a structure of a closed symmetric
monoidal category on the category of cochain complexes over $\cat C$. 

Let $X^*$ and $Y^*$ be two cochain complexes over $\cat C$. Let $X^*
\otimes Y^*$ be the cochain complex with $(X^* \otimes Y^*)^n =
\bigoplus_{p + q = n} X^p \otimes Y^q$ and such that the differential
$$d^n: \bigoplus_{p + q = n} X^p \otimes Y^q \to \bigoplus_{p + q = n
  + 1} X^p \otimes Y^q$$
is given componentwise by $$\diff^p \otimes
\id_{Y^q} + (-1)^p \id_{X^p} \otimes \diff^q: X^p \otimes Y^q \to X^{p
  + 1} \otimes Y^q \oplus X^p \otimes Y^{q + 1}.$$
On the defines
tensor products of morphisms between cochain complexes over $\cat C$
in the obvious way. By defining associators in the most obvious way,
this makes $\cdot \otimes \cdot: \cat C^* \times \cat C^* \to \cat
C^*$ a tensor product in the category of cochain complexes over $\cat
C$.

A unit for this tensor product is given by the cochain
complex $\unit^*$ with $\unit^n = 0$ for $n \neq 0$ and $\unit^0 =
\unit$. The unit laws are defined in the most obvious way. 

Let $\gamma^*: X^* \otimes Y^* \to Y^* \otimes X^*$ be the isomorphism
natural in $X^*$ and $Y^*$ that is given componentwise given by
$(-1)^{pq} \gamma: X^p \otimes Y^q \to Y^q \otimes X^p$ where $\gamma: X^p
\otimes Y^q \to Y^q \otimes X^p$ is given by the braiding of $\cat
C$. This gives a symmetric braiding for the category $\cat C^*$. 

Finally, an inner hom $\hom^*(Y^*, Z^*)$ for the tensor product on
$\cat C^*$ is defined by $\hom^n(Y^*, Z^*) = \prod_{p + n = q}
\hom(X^p, Y^q)$ and such that the differential $$\diff^n: \prod_{p + n
  = q} \hom(X^p, Y^q) \to \prod_{p + n + 1 = q} \hom(X^p, Y^q)$$
is
given componentwise by
\begin{multline*}
  \hom(X^p, \diff^q) - (-1)^n \hom(\diff^{p -
    1}, Y^q): \hom(X^p, Y^q) \\
  \longrightarrow \hom(X^p, Y^{q + 1}) \oplus \hom(X^{p - 1}, Y^q).
\end{multline*}
The adjunction morphisms are defined in the most obvious way.

\begin{proposition}
  With the so defined structures, the category of cochain complexes
  over a closed symmetric monoidal additive category becomes naturally
  a closed symmetric monoidal additive category. 
\end{proposition}
\qed

\subsection{Monoids}

Let $\cat C$ be a monoidal category. 

Let $A$ be an object of $\cat C$. A \emph{multiplication law for $A$}
is a morphism
\[\mu: A \otimes A \to A.\]
It is \emph{associative} if the canonical diagram
\[\xymatrix{
  A \otimes A \otimes A \ar[r]^{\id \otimes \mu} \ar[d]_{\mu \otimes
    \id} &
  A \otimes A \ar[d]^\mu \\
  A \otimes A \ar[r]_\mu & A
}\]
commutes. 

A morphism $\eta: \unit \to A$ is a \emph{unit for the multiplication
  law} if the canonical diagrams
\[\xymatrix{
  A \ar[r]^{\id \otimes \eta} \ar@{=}[dr] & A \otimes A \ar[d]^\mu \\
  & A
}\]
and
\[\xymatrix{
  A \ar[r]^{\eta \otimes \id} \ar@{=}[dr] & A \otimes A \ar[d]^\mu \\
  & A 
}\]
commute. 

\begin{definition}
  A \emph{monoid in $\cat C$} or a \emph{unital associative algebra in
    $\cat C$} is an object in $\cat C$ together with an associative
  multiplication law and a unit for this law.
\end{definition}

Assume that $\cat C$ is symmetric as a monoidal category. Let $\gamma:
A \otimes A \to A \otimes A$ be the canonical morphism that
interchanges the two factors $A$. 

The multiplication law $\mu$ is \emph{commutative} if the canonical
diagram
\[\xymatrix{
  A \otimes A \ar[r]^\gamma \ar[dr]_\mu & A \otimes A \ar[d]^\mu \\
  & A
}\]
commutes. 

\begin{definition}
  A \emph{commutative monoid in $\cat C$} or a \emph{unital
    commutative algebra in $\cat C$} is an object in $\cat C$ together
  with an associative commutative multiplication law and a unit for
  this law.
\end{definition}

Let $R$ be a unital (commutative) ring. 
\begin{example}
  A unital (commutative) algebra in the category of left $R$-modules
  is nothing else than an ordinary (commutative) $R$-algebra. 
\end{example}

\begin{remark}
  We leave the obvious definition of a morphism of (commutative)
  monoids in $\cat C$ and the composition law for these morphisms to
  the reader. In particular, one can define the category of
  (commutative) monoids in $\cat C$.
\end{remark}

\subsection{Modules over monoids}

Let $\cat C$ be a monoidal category and let $A$ be a monoid in $\cat
C$ with the multiplication law $\mu: A \otimes A \to A$. 

Let $M$ be any object. An \emph{$A$-operation on $M$} is given by a
morphism
\[\nu: A \otimes M \to M.\]

This operation is \emph{associative} if the canonical diagram
\[\xymatrix{
  A \otimes A \otimes M \ar[r]^{\id \otimes \nu} \ar[d]_{\mu \otimes
    \id} &
    A \otimes M \ar[d]^{\nu} \\
    A \otimes M \ar[r]_{\nu} & M
}\]
commutes. 

It is \emph{unital} if the canonical diagram
\[\xymatrix{
  M \ar[r]\ar@{=}[dr] & A \otimes M \ar[d] \\
  & M
}\]
commutes. 

\begin{definition}
  A \emph{(left) $A$-module} or a \emph{(left) module over $A$} is an
  object of $\cat C$ together with a unital and associative
  $A$-operation. 
\end{definition}

Let $R$ be a unital ring and let $A$ be an $R$-algebra, which we also
view as a unital associative algebra in the category of left
$R$-modules. 
\begin{example}
  A left module over the monoid $A$ is nothing else than an ordinary
  left $A$-module. 
\end{example}

\begin{remark}
  We leave the obvious definition of a morphism of left modules over a
  monoid and the composition law for these morphisms to the reader. In
  particular, one can define the category of left modules over an
  algebra.
\end{remark}

\section{Operads, Algebras and Modules}

\subsection{Permutations}

We use the latter $\SG$ to denote the permutation groups:
\begin{definition}
  The \emph{permutation group in $n$ letters} is the group of
  bijections of the set $\{1, \dots, n\}$ into itself and denoted by
  $\SG_n$. 
\end{definition}

Let $m_1, \dots, m_n$ be non-negative integers and set $m := \sum_{i =
  1}^n m_i$. Let $\tau_i \in \SG_{m_i}$ be permutations. These induce
a permutation $\tau \in \SG_m$ which is given by
\[\tau\left(\sum_{i = 1}^{j - 1} m_j + k\right) = \tau_j(k) + \sum_{i = 1}^{j -
  1} m_j
\]
for $1 \leq k \leq m_j$. 
\begin{definition}
  The permutation $\tau$ defined above is the \emph{sum of the
    permutations $\tau_1, \dots, \tau_n$} and denoted by $\tau_1 +
  \dots + \tau_n$.
\end{definition}

Let $\sigma \in \SG_n$ be a permutation. It induces a permutation
$\tilde \sigma \in \SG_m$ which is given by
\[\tilde\sigma\left(\sum_{i = 1}^{j - 1} m_j + k\right) =
\sum_{\sigma(i) < \sigma(j)} m_i + k\]
for $1 \leq k \leq m_j$. 
\begin{definition}
  The permutation $\tilde \sigma$ defined above is the \emph{block
    permutation induced by $\sigma$ for the blocks given by $m_1,
    \dots, m_n$} and denoted by $\sigma_{m_1, \dots, m_n}$. 
\end{definition}

\subsection{Operads}
Operads describe certain types of algebras. There is, e.g., the
\emph{associative operad} that governs associative algebras. Another
example is the \emph{commutative operad} governing commutative
algebras. 

For what follows, let $\cat C$ be a closed symmetric monoidal
category, e.g.~the category of $k$-modules for a commutative ring with
unit $k$. 

\begin{definition}
  An \emph{$\SG$-module $\sgm M$} is a sequence $(\sgm M(n))_{n \in
    \mathbf N_0}$ of objects $\sgm M(n)$ in $\cat C$ together with a
  right action of $\SG_n$ on each $\sgm M(n)$.
\end{definition}

Let $\sgm M$ be an $\SG$-module. An \emph{operadic composition law} on
$\sgm M$ is given by a family of morphisms
\[\gamma: \sgm M(n) \otimes \sgm M(m_1) \otimes \dots \otimes \sgm
M(m_n) \to M(m)\] with $m := \sum_{i = 1}^n m_i$. 

Such an operadic composition law is called \emph{associative} if the
canonical diagrams
\[\xymatrix{
  \sgm M(n) \otimes \bigotimes_{i = 1}^n \left(\sgm M(m_i) \otimes
  \bigotimes_{j = 1}^{m_i} \sgm M(l_{i, j})\right)
  \ar[r] \ar[d] &
  \sgm M(n) \otimes \bigotimes_{i = 1}^n \sgm M(l_i) \ar[d] \\
  \sgm M(m) \otimes \bigotimes_{i = 1}^n \bigotimes_{j = 1}^{m_i} \sgm
  M(l_{i, j}) \ar[r] & \sgm M(l)
}\]
with $l_i := \sum_{j = 1}^{m_i} l_{i,j}$ and $l := \sum_{i = 1}^n l_i$
built up from the composition law commute. 

The operadic composition law is called \emph{equivariant} if the
canonical diagrams
\[\xymatrix@C+6pt{
  \sgm M(n) \otimes \left(\sgm M(m_1) \otimes \dots \otimes \sgm
    M(m_n)\right)
  \ar[r]^{(\cdot \sigma) \otimes (\sigma^{-1} \cdot)} \ar[d] &
  \sgm M(n) \otimes \left(\sgm M(m_{\sigma(1)}) \otimes \dots
    \otimes \sgm M(m_{\sigma_n})\right)
  \ar[d] \\
  \sgm M(m) \ar[r]_{\sigma_{m_{\sigma(1)}, \dots, m_{\sigma(n)}}} &
  \sgm M(m)
}\]
for all $\sigma \in \SG_n$
and
\[\xymatrix@C+24pt{
  \sgm M(n) \otimes \left(\sgm M(m_1) \otimes \dots \otimes \sgm
    M(m_n)\right)
  \ar[r]^{\id \otimes (\cdot \tau_1) \otimes \dots \otimes (\cdot
    \tau_n)} \ar[d] & \sgm M(n) \otimes \left(\sgm M(m_1) \otimes \dots
    \otimes \sgm M(m_n)\right)
  \ar[d] \\
  \sgm M(m) \ar[r]_{\cdot (\tau_1 + \dots + \tau_n)} & \sgm M(m)
}\]
for all $\tau_i \in \SG_{m_i}$ commute.

A morphism $\eta: \unit \to \sgm M(1)$ is a \emph{unit for the
  operadic composition law} if the canonical diagrams
\[\xymatrix@C+24pt{
  \sgm M(n) \ar[r]^{\id \otimes \eta^{\otimes n}} \ar@{=}[dr] &
  \sgm M(n) \otimes \sgm M(1)^{\otimes n} \ar[d] \\
  & \sgm M(n)
}\]
and
\[\xymatrix@C+24pt{
  \sgm M(m) \ar[r]^{\eta \otimes \id} \ar@{=}[dr] &
  \sgm M(1) \otimes \sgm M(m) \ar[d] \\
  & \sgm M(n)
}\]
commute. 

We give the definition of an operad in a second. Before that, however,
we want to note that the reader who wants to learn more about operads
should consult the monograph~\cite{MarklShniderStasheff02} and the
references therein.
\begin{definition}
  An \emph{operad} is an $\SG$-module together with an associative and
  equivariant operadic composition law and a unit for this law.
\end{definition}

The prototype of an operad is the so called \emph{endomorphism operad
  $\End_V$ over an object $V$ in $\cat C$}: Its underlying
$\SG$-module is given by $\End_V(n) = \hom(V^{\otimes n}, V)$ with the
canonical operations of each $\SG_n$ on the right. There is an
equivariant operadic composition law on $\End_V$ given by the ordinary
composition morphisms
$$\End_V(n) \otimes \End_V(m_1) \otimes \dots \otimes \End_V(m_n) \to
\End_V(m)$$
with $m := \sum_{i = 1}^n m_i$.
  
A canoncial unit for the operadic composition law on $\End_V$ is given
by the morphism $\unit \to \End_V(1) = \hom(V, V)$ which corresponds
to the identity morphism on $V$.

\begin{example}
  The operad $\End_V$ given by the $\SG$-module and the operadic
  composition law defined above is the \emph{endomorphism operad over
    $V$}.
\end{example}

Consider the $\SG$-module $\Ass$ which is given by $\Ass(n) =
\unit^{\amalg \SG_n}$ with the canonical operations of each $\SG_n$ on
the right. The natural map from $\unit \to \Ass(n)$ corresponding to
the inclusion of the summand indexed by a permutation $\sigma \in
\SG_n$ is denoted by $\cdot \sigma: \unit \to \Ass(n)$.

There is a unique equivariant operadic composition law $\gamma$ on $\Ass$ such
that the canonical diagrams
\[\xymatrix@C+48pt{
  \unit \otimes \unit^{\otimes n} \ar[r]\ar[d]_{(\cdot \id) \otimes
    (\cdot \id)^n} &
  \unit \ar[d]^{\cdot \id} \\
  \Ass(n) \otimes \Ass(m_1) \otimes \dots \otimes \Ass(m_n)
  \ar[r]_{\gamma} & \Ass(m) \\
}\] commute.

A canonical unit for the operadic composition law on $\Ass$ is given
by $\cdot \id: \unit \to \Ass(1)$. 

\begin{example}
  The operad $\Ass$ given by the $\SG$-module and the operadic
  composition law defined above is the \emph{associative operad}.
\end{example}

Consider the $\SG_n$-module $\Com$ which is given by $\Com(n) = \unit$
with the trivial operation of each $\SG_n$ on the right. A canonical
operadic composition law is given by the unit laws of the underlying
symmetric monoidal category. This operadic composition law is
trivially equivariant and has a canonical unit.

\begin{example}
  The operad $\Com$ given by the $\SG$-module and the operadic
  composition law defined above is the \emph{commutative operad}.
\end{example}

\begin{remark}
  We leave the obvious definition of a morphism of operads and the
  composition law for these morphisms to the reader.  In particular,
  one can define the category of operads in $\cat C$.
\end{remark}

\subsection{Algebras}

In this subsection, we define the notion of an algebra over an
operad. 

Again let $\cat C$ denote a closed symmetric monoidal category. Let
$\sgm O$ be an operad in $\cat C$. 

Let $A$ be an object in $\cat C$. An \emph{$\sgm O$-multiplication law
  on $A$} is a family of morphisms
\[
\sgm O(n) \otimes A^{\otimes n} \to A. 
\]

Such a multiplication law is called \emph{associative} if the
canonical diagrams
\[\xymatrix{
  \sgm O(n) \otimes \sgm O(m_1) \otimes \dots \otimes \sgm O(m_n)
  \otimes A^{\otimes m} \ar[r]\ar[d] &
  \sgm O(n) \otimes A^{\otimes n} \ar[d] \\
  \sgm O(m) \otimes A^{\otimes m} \ar[r] & A }\] built up from the
operadic composition law and the $\sgm O$-multiplication law with $m
:= \sum_{i = 1}^n m_i$ commute.

The multiplication law is called \emph{equivariant} if the canonical
diagrams
\[\xymatrix@C+24pt{
  O(n) \otimes A^{\otimes n}
  \ar[r]^{(\cdot \sigma) \otimes (\sigma^{-1} \cdot)} \ar[dr] &
  O(n) \otimes A^{\otimes n} \ar[d] \\
  & A
}\]
for all $\sigma \in \SG_n$ commute. 

The multiplication law is called \emph{unital} if the canonical
diagram
\[\xymatrix{
A \ar[r] \ar@{=}[dr] & \sgm O(1) \otimes A \ar[d] \\
& A
}\]
commutes. 

\begin{definition}
  An \emph{algebra over $\sgm O$} is an object of $\cat C$ together
  with a unital, equivariant, and associative $\sgm O$-multiplication
  law. 
\end{definition}
(For more on algebras over operads, one may
consult~\cite{MarklShniderStasheff02}.)

\begin{remark}
  In this monograph, one can also find the one-to-one correspondence
  between algebra structures over $\sgm O$ on an object $A$ on the one
  hand and morphisms $\sgm O \to \End_A$ of operads in $\sgm C$ on the
  other hand.
\end{remark}

Let $A$ be a unital, associative algebra in $\cat C$. There is a
unique equivariant $\Ass$-multiplication law $\gamma$ on $A$ such that
the canonical diagrams
\[\xymatrix@C+24pt{
  A^{\otimes n} \ar[r]^{(\cdot \id) \otimes \id^{\otimes n}} \ar[dr] &
  \Ass(n) \otimes A^{\otimes n} \ar[d]^{\gamma} \\
  & A }\] commute. This law is associative and unital. Thus we get the
following:

\begin{example}
  Each unital, associative algebra in $\cat C$ defines canonically an
  algebra over $\Ass$. 

  (In fact, this defines a one-to-one correspondence between unital,
  associative algebras and algebras over $\Ass$.) 
\end{example}

Let $C$ be a unital, commutative algebra in $\cat C$. There exists a
unique $\Com$-multiplication law $\gamma$ on $C$ such that the canonical
diagrams
\[\xymatrix{
  C^{\otimes n} \ar[r] \ar[dr] &
  \Com(n) \otimes C^{\otimes n} \ar[d]^{\gamma} \\
  & C
}\]

commute. This law is associative, equivariant, and unital. Thus we get
the following:
\begin{example}
  Each unital, commutative algebra in $\cat C$ defines canonically an
  algebra over $\Com$. 

  (In fact, this defines a one-to-one correspondence between unital,
  commutative algebras and algebras over $\Com$.) 
\end{example}

\begin{remark}
  Again, we leave the obvious definition of a morphism of algebras
  over an operad and the composition law for these morphisms to the
  reader. In particular, one can define the category of algebras over
  an operad.
\end{remark}

\subsection{The free algebra}

There is a ``free construction'' in the category of algebras over a
fixed operad $\sgm O$ in a cocomplete closed symmetric monoidal
category $\cat C$. 

Let $V$ be an object in $\cat C$. Each $\sigma \in \SG_n$ induces a
morphism $$(\cdot \sigma) \otimes (\sigma^{-1} \cdot): \sgm O(n)
\otimes V^{\otimes n} \to \sgm O(n) \otimes V^{\otimes n}.$$
This
defines a canonical operation of $\SG_n$ on $\sgm O(n) \otimes
V^{\otimes n}$. The corresponding object of coinvariants is denoted by
$\sgm O(n) \otimes_{\SG_n} V^{\otimes n}$.

We set
\[
F_{\sgm O} V := \coprod_{n = 0}^\infty \sgm O(n) \otimes_{\SG_n}
V^{\otimes n}. 
\]
There is a unique $\sgm O$-multiplication law on $F_{\sgm O} V$ which
is induced by the morphisms given componentwise by
\[
\sgm O(n) \otimes (\sgm O(m_1) \otimes V^{\otimes m_1}) \otimes \dots
\otimes (\sgm O(m_n) \otimes V^{\otimes m_n}) \longrightarrow \sgm
O(m) \otimes V^{\otimes m}\] with $m := \sum_{i = 1}^n m_i$ that are
induced by the operadic composition law on $\sgm O$. This
multiplication law is associative, equivariant, and unital.
\begin{definition}
  The \emph{free $\sgm O$-algebra over $V$} is
  the $\sgm O$-algebra $F_{\sgm O} V$ given by the $\sgm
  O$-multiplication law defined above. 
\end{definition}
(More on this free construction can be found
in~\cite{MarklShniderStasheff02}.)

There is a natural morphism $V \to F_{\sgm O} V$ that is induced by
the unit $\unit \to \sgm O(1)$ of the operad $\sgm O$. 

\begin{proposition}
  The free algebra over $V$ is in fact \emph{free}, i.e.~it is the
  initial object in the comma category $(V, \#)$ where $\#$ is the
  forgetful functor from the category of $\sgm O$-algebras to $\cat
  C$.
\end{proposition}
\qed

\begin{example}
  The free $\Ass$-algebra over an object is naturally isomorphic to
  the tensor algebra over that object. 
\end{example}

\begin{example}
  The free $\Com$-algebra over an object is naturally isomorphic to
  the symmetric algebra over that object. 
\end{example}

\subsection{Modules}

Let $\sgm O$ be an operad over a closed symmetric monoidal category
$\cat C$. Let $A$ be an $\sgm O$-algebra. 

Let $M$ be an object in $\cat C$. An $A$-operation on $M$ is a family
of morphisms
\[\sgm O(n + 1) \otimes A^{\otimes n} \otimes M \to M.\]

Such an operation is \emph{associative} if the canonical diagrams
\[\xymatrix{
  \sgm O(n + 1) \otimes \sgm O(m_1) \otimes \dots
  \otimes \sgm O(m_{n + 1})
  \otimes A^{\otimes m} \otimes M \ar[r]\ar[d] &
  \sgm O(n + 1) \otimes A^{\otimes n} \otimes M \ar[d] \\
  \sgm O(m + 1) \otimes A^{\otimes m} \otimes M \ar[r] & M
}\]
with $m := \sum_{i = 1}^{n + 1} m_i - 1$ commute. 

The operation is \emph{equivariant} if the canonical diagrams
\[\xymatrix@C+48pt{
  O(n + 1) \otimes A^{\otimes n} \otimes M
  \ar[r]^{(\cdot (\sigma + \id)) \otimes (\sigma^{-1} \cdot)
    \otimes \id} \ar[dr] &
  O(n + 1) \otimes A^{\otimes n} \otimes M\ar[d] \\
  & M
}\]
for all $\sigma \in \SG_n$ commute. 

The operation is \emph{unital} if the canonical diagram
\[\xymatrix{
  M \ar[dr] \ar[r] & \sgm O(1) \otimes M \ar[d] \\
  & M
}\]
commutes. 

\begin{definition}
  A \emph{module over $A$ (or an $A$-module)} is an object of $\cat C$
  together with a unital, equivariant, and associative $A$-operation
  on it. 
\end{definition}

\begin{example}
  The algebra $A$ is naturally a module over itself. 
\end{example}

Let $A$ be a unital, associative algebra in $\cat C$, which we also
view as an $\Ass$-algebra, and let $M$ be an $A$-\emph{bi}module.
There is a unique equivariant $A$-operation $\gamma$ on $M$ such that
the canonical diagrams
\[\xymatrix@C+48pt{
  A^{\otimes n} \otimes M \ar[r]^{(\cdot \id) \otimes \id^{\otimes n}
    \otimes \id}  \ar[dr] &
  \Ass(n + 1) \otimes A^{\otimes n} \otimes M \ar[d]^{\gamma} \\
  & M
}\]
commute. This law is associative and unital. This can be summarised as
follows:
\begin{example}
  Each bimodule over $A$ is canonically an $A$-module when $A$ is
  viewed as an $\Ass$-algebra.
\end{example}

\begin{example}
  Similarly, one can show that $C$-modules for a unital, commutative
  algebra $C$ are canonically modules over $C$ when $C$ is viewed as a
  $\Com$-algebra.
\end{example}

\begin{remark}
  Again, we leave the obvious definition of a morphism of modules over
  an algebra and the composition law for these morphisms to the
  reader. In particular, one can define the category of modules over
  an algebra.
\end{remark}

\begin{remark}
  Assume that $\cat C$ is a closed symmetric monoidal additive
  category. This defines a canonical structure of an additive category
  on the category of modules over an algebra over an operad in $\cat
  C$. In case $\cat C$ is abelian, this structure is also abelian. One
  way to show this is to show that the category of modules is
  isomorphic to a category of left modules over a certain associative
  algebra $U$ (which is defined in the next subsection).
\end{remark}

\subsection{The universal enveloping algebra}

In this subsection, we recall the notion of the universal enveloping
algebra, which is an ordinary unital, associative algebra, of an
algebra over some operad. For further reading, we
suggest~\cite{Hinich97} and the references therein. Let $\cat C$ be a
cocomplete closed symmetric monoidal category and $\sgm O$ an operad
in $\cat C$. Let $A$ be an $\sgm O$-algebra.

Set
\[T(A) := \coprod_{n = 0}^\infty \sgm O(n + 1) \otimes_{\SG_n}
A^{\otimes n},\]
which is called the \emph{tensor algebra over $A$}.

Let $U(A)$ be the maximal quotient of $T(A)$ such that the canonical
diagrams
\[\xymatrix{
  \left(\sgm O(n + 1) \otimes \sgm O(m_1) \otimes \dots \otimes \sgm
    O(m_{n + 1})\right) \otimes_{\SG_m} A^{\otimes m} \ar[r]\ar[d] &
  \sgm O(n + 1) \otimes_{\SG_n} A^{\otimes n} \ar[d] \\
  \sgm O(m + 1) \otimes_{\SG_m} A^{\otimes m} \ar[r] & U(A) }\] with
$m := \sum_{i = 1}^{n + 1} m_i - 1$ commute.

There is a unique structure of an associative algebra on $T(A)$ which
is induced by the canonical morphisms given componentwise by
\begin{multline*}
  (\sgm O(n + 1) \otimes A^{\otimes n}) \otimes (\sgm O(m + 1) \otimes
  A^{\otimes m}) \longrightarrow \sgm O(n + 1) \otimes \sgm
  O(1)^{\otimes n} \otimes \sgm O(m + 1) \otimes A^{\otimes (n + m)}
  \\
  \longrightarrow \sgm O(n + m + 1) \otimes A^{\otimes (n + m)}
\end{multline*}
The canonical morphism $\unit \to T(A)$ makes up a canonical
unit for $T(A)$. This structure makes $U(A)$ naturally into a
unital, associative algebra in $\cat C$ such that the quotient morphism
$T(A) \to U(A)$ is a morphism of unital, associative algebras. 

\begin{definition}
  The unital associative algebra $U(A)$ in $\cat C$ as defined above
  is the \emph{universal enveloping algebra of $A$}. 
\end{definition}

\begin{example}
  Let $A$ be a classical unital associative algebra in $\cat C$,
  which we also view as an $\Ass$-algebra. The universal enveloping
  algebra over $A$ is given by $A \otimes A^\op$. 
\end{example}

\begin{example}
  Let $A$ be a classical unital commutative algebra in $\cat C$,
  which we also view as an $\Com$-algebra. The universal envelopping
  algebra over $A$ is given by $A$ itself. 
\end{example}

\begin{remark}
  As already stated in the previous subsection, the category of
  modules over the associative algebra $U(A)$ is isomorphic to the
  category of $A$-modules. In fact, an $A$-module $E$ is the same as a
  morphism $U(A) \to \hom(E, E)$ of associative algebras in $\cat C$.
\end{remark}

Let $W$ be an object in $\cat C$. Set
\[F_A(W) := U(A) \otimes W.\] There is a unique equivariant
$A$-operation on $F_A(W)$ which is induced by the natural morphisms
given componentwise by
\begin{multline*}
  \sgm O(n + 1) \otimes A^{\otimes n} \otimes (\sgm O(m + 1) \otimes
  A^{\otimes m} \otimes W)
  \\
  \longrightarrow \sgm O(n + 1) \otimes \sgm O(1)^{\otimes
    n} \otimes \sgm O(m + 1) \otimes A^{\otimes (n + m)} \otimes W
  \\
  \longrightarrow 
  \sgm O(n + m + 1) \otimes A^{\otimes (n + m)} \otimes W. 
\end{multline*}

This $A$-operation is associative and unital. 

\begin{definition}
  The \emph{free $A$-module over $W$} is the $A$-algebra $F_A(W)$
  given by the $A$-operation defined above. 
\end{definition}
(See also~\cite{Hinich97}.)

There is a natural morphism $W \to F_A(W)$ induced by the canonical
morphism $\unit \to U(A)$ (the unit) induced by the unit $\unit \to
\sgm O(1)$ of the operad $\sgm O$.

\begin{proposition}
  The free $A$-module over $W$ is in fact \emph{free}, i.e.~it is the
  initial object in the comma category $(W, \#)$, where $\#$ is the
  forgetful functor from the category of $A$-modules to $\cat C$. 
\end{proposition}
\qed

\subsection{Derivations}

Let $\cat C$ be a closed symmetric monoidal additive category and $\sgm O$ an
operad in $\cat C$. Let $A$ be an $\sgm O$-algebra. 

Let $M$ be an $A$-module. A morphism $d: A \to M$ is \emph{derivative}
if the canonical diagrams
\[\xymatrix@C+72pt{
  \sgm O(n) \otimes A^n \ar[r]^{\id \otimes \sum_{p + 1 + q = n}
    \id^{\otimes p} \otimes d \otimes \id^{\otimes q}} \ar[d] &
  \sgm O(n) \otimes A^n \ar[d] \\
  A \ar[r]_d & M }\] commute.

\begin{definition}
  A \emph{derivation of $A$ (into the module $M$)} is a derivative
  morphism as defined above.
\end{definition}

\begin{example}
  Let $A$ be a unital associative algebra in $\cat C$, which we also
  view as an $\Ass$-algebra and let $M$ be an $A$-module. Every
  derivation $d: A \to M$ in the classical sense is a derivation in
  the above sense. 
\end{example}

\begin{example}
  An analogous result holds true in the category of commutative
  algebras in $\cat C$, i.e.~in the category of algebras over $\Com$. 
\end{example}

Let $V$ be an object in $\cat C$ and $M$ an $A$-module. Let $\phi: V
\to M$ be a morphism in $\cat C$. There is exactly one derivation
$d_\phi: F_{\sgm O}(V) \to M$ such that the following canonical
diagram
\[\xymatrix{
  V \ar[r] \ar[dr]_{\phi} &
  F_{\sgm O} V \ar[d]^{d_\phi} \\
  & M
}\]
commutes. 

\begin{proposition}
  The mapping $\phi \mapsto d_\phi$ sets up a natural bijection
  between the morphisms in $\cat C$ from $V$ to $M$ and the
  derivations from $F_{\sgm O}(V)$ into $M$. 
\end{proposition}
\qed

\begin{example}
  Let $V$ and $W$ be two objects in $\cat C$. Set $A := F_{\sgm
    O}(V)$. Every morphism from $V$ to $W$ in $\cat C$ induces a
  natural derivation $A \to F_A(W)$. 
\end{example}

Set $\tilde d := \eta \otimes \id: A \to U(A) \otimes A$, where $\eta:
\unit \to U(A)$ is the canonical unit of the universal envelopping
algebra. Let $\Omega_A$ be the maximal $A$-module quotient of $F_A(A)
= U(A) \otimes A$ such that the canonical diagrams
\[\xymatrix{
  & \sgm O(n) \otimes A^{\otimes n} \ar[dl] \ar[dr]^{\hspace{24pt}\id
    \otimes \sum_{p + 1 + q = n} \id^{\otimes p}
    \otimes \tilde d \otimes \id^{\otimes q}} \\
  A \ar[d]_{\tilde d} & & \sgm O(n) \otimes A^{\otimes (n - 1)}
  \otimes (U(A) \otimes A)
  \ar[d] \\
  U(A) \otimes A \ar[dr] & & U(A) \otimes A \ar[dl] \\
  & \Omega_A }\] commute.

Let $d: A \to \Omega_A$ be the composition of $\tilde d$ and the
quotient map $U(A) \otimes A \to \Omega_A$. The morphism $d$ is
derivative. 

\begin{definition}
  The derivation $d: A \to \Omega_A$ is the \emph{K\"ahler (or
    universal) derivation} of $A$. 
\end{definition}

\begin{proposition}
  The universal derivation $d: A \to \Omega_A$ is in fact
  \emph{universal}, i.e.~it is the initial object in the category of
  derivations of $A$. 
\end{proposition}
\qed

\begin{example}
  Let $A$ be a unital, associative algebra in $\cat C$, which we also
  view as an $\Ass$-algebra. The universal derivation is the classical
  one: set $$\tilde d := \eta \otimes \eta^\op \otimes \id: A \to A
  \otimes A^\op \otimes A,$$ where $\eta$ and $\eta^\op$ are the units
  of $A$ and $A^\op$, respectively. Then $\Omega_A$ is the maximal
  quotient of $A \otimes A^\op \otimes A$ such that the canonical
  diagram
  \[\xymatrix{
    & A \otimes A
    \ar[dl] \ar[dr]^{\hspace{12pt}\tilde d \otimes \id + \id \otimes \tilde d} \\
    A \ar[d]_{\tilde d} & & A \otimes (A \otimes A^\op \otimes A)
    \ar[d] \\
    A \otimes A^\op \otimes A \ar[dr] & & A \otimes A^\op \otimes A \ar[dl] \\
    & \Omega_A }\] commutes. The universal derivation $d: A \to
  \Omega_A$ is given by the composition of $\tilde d$ and the quotient
  map $A \otimes A^\op \otimes A \to \Omega_A$.
\end{example}

\begin{example}
  In the case of unital commutative algebras, the classical and the
  operadic notion of the K\"ahler derivation coincide as well. 
\end{example}

\begin{example}
  Let $V$ be an object of $\cat C$. Let $A := F_{\sgm O}(V)$ be the
  free $\sgm O$-algebra over $V$. The natural morphism $d: A \to
  U(A) \otimes V$ is the universal derivation of $A$. 
\end{example}

\subsection{Lax products}

Given modules over an operad, one may consider their tensor product in
the underlying closed symmetric monoidal category. In the case of the,
say, commutative operad, this product, however, does not coincide with
the tensor product as modules. Thus we propose here a new product of
modules over a general algebra that is closer to the given module
structure than the tensor product in $\cat C$.

Let $\cat C$ be a cocomplete closed symmetric monoidal category and
let $\sgm O$ be an operad in $\cat C$. Let $A$ be an $\sgm O$-algebra. 

Let $M_1, \dots, M_m$ be modules over $A$. Let $P_A(M_1, \dots, M_m)$
be the maximal quotient of $\coprod_{n = 0}^\infty \sgm O(n + m)
\otimes_{\SG_n} A^{\otimes n} \otimes M_1 \otimes \dots \otimes M_m$
such that the canonical diagrams
\[\xymatrix{
  {\begin{aligned}
      \left(\sgm O(n + m) \otimes \sgm O(k_1) \otimes \dots \otimes \sgm O(k_{n
          + m})\right) \\
      \otimes_{\SG_k} A^{\otimes k}
      \otimes M_1 \otimes \dots
      \otimes M_m
    \end{aligned}}
  \ar[r]\ar[d] &
  {\begin{aligned}
      \sgm O(n + m)
      \otimes_{\SG_n} A^{\otimes n} \\ \otimes M_1 \otimes
      \dots \otimes M_m
    \end{aligned}}
  \ar[d] \\
  \sgm O(k + m) \otimes_{\SG_k} A^{\otimes k} \otimes M_1 \otimes
  \dots \otimes M_m \ar[r] &
  P_A(M_1, \dots, M_m)
}\]
with $k := \sum_{i = 1}^{n + m} k_i - m$ commute. 

There is a natural (unital, associative) $A$-operation on $P_A(M_1,
\dots, M_n)$ which is induced by the canonical morphisms given
componentwise by
\begin{multline*}
  \sgm O(k + 1) \otimes A^{\otimes k} \otimes (\sgm O(n + m) \otimes
  A^{\otimes n} \otimes M_1 \otimes \dots \otimes M_m)
  \\
  \longrightarrow \sgm O(k + 1) \otimes (\sgm O(1))^{\otimes k}
  \otimes \sgm O(n + m) \otimes A^{\otimes (k + n)} \otimes M_1
  \otimes \dots \otimes M_m
  \\
  \longrightarrow \sgm O(k + n + m) \otimes A^{\otimes (k + n)}
  \otimes M_1 \otimes \dots \otimes M_m.
\end{multline*}

\begin{definition}
  The $A$-module $P_A(M_1, \dots, M_m)$ as defined above is the
  \emph{lax product of the $A$-modules $M_1, \dots, M_m$}. 
\end{definition}

\begin{example}
  Let $M$ be an $A$-module. It is $P_A() = A$ and $P_A(M) = M$ as
  $A$-modules. 
\end{example}

\begin{remark}
  There are natural morphisms
  \[P_A(P_A(M_{1,1}, \dots, M_{1, m_1}), \dots, P_A(M_{n, 1}, \dots, M_{n,
    m_n})) \longrightarrow P_A(M_{1,1}, \dots, M_{n, m_n}),
  \]
  and
  \[P_A(M_1, \dotsm M_n) \longrightarrow P_A(M_{\sigma(1)}, \dots, M_{\sigma(n)}),\]
  $\sigma \in \SG_n$, that make certain coherence diagrams
  commutative.
  
  (In fact, the category of $A$-modules with the lax product becomes a
  \emph{lax symmetric monoidal category}, for this notion see,
  e.g.,~\cite{DayStreet03}.)
\end{remark}

\begin{definition}
  Let $M$ and $E_1, \dots, E_m$ be $A$-modules. We set
  \[P_A^n(M, E_1, \dots, E_m) := P_A\left(\underbrace{M, \dots,
      M}_{\textrm{$n$-times}}, E_1, \dots, E_m\right).\] The object of
  $\SG_n$-coinvariants by the natural $\SG_n$-operation on $P_A^n(M,
  E_1, \dots, E_m)$ is denoted by $S_A^n(M, E_1, \dots, E_m)$.
\end{definition}

\begin{example}
  Let $A$ be a unital, associative algebra in $\cat C$, which we also
  view as an $\Ass$-algebra. Let $M_1, \dots, M_n$ be modules over
  $A$. Their lax product is given by
  \[P_A(M_1, \dots, M_n) = \coprod_{\sigma \in \SG_n} M_{\sigma(1)}
  \otimes_A \dots \otimes_A M_{\sigma(n)}\]
  with the obvious structure of a module over $A$. 
\end{example}

\begin{example}
  Let $C$ be a unital, associative algebra in $\cat C$, which we also
  view as a $\Com$-algebra. Let $M_1, \dots, M_n$ be modules over
  $C$. Their lax product is given by
  \[P_C(M_1, \dots, M_n) = M_1 \otimes_C \dots \otimes_C M_n\] with the
  obvious structure of a module over $C$. 
\end{example}

\begin{example}
  Let $W_1, \dots, W_m$ be objects in $\cat C$. The lax tensor product
  of the free $A$-modules $M_i := F_A(W_i) = U(A) \otimes W_i$ is
  given by
  \[P_A(M_1, \dots, M_m) = U(A) \otimes W_1 \otimes \dots \otimes W_m
  = F_A(W_1 \otimes \dots \otimes W_m).\] 
\end{example}

\subsection{Lax inner hom's}

In this subsection, we want to construct right adjoints to the various
lax products defined in the previous subsection. 

Let $\cat C$ be a bicomplete closed symmetric monoidal category, and
let $\sgm O$ be an operad in $\cat C$. Let $A$ be an $\sgm O$-algebra. 

Let $M_2, \dots, M_m$ and $N$ be modules over $A$. Let $H_A(M_2, \dots,
M_m; N)$ be the maximal subobject of $\prod_{n \ge 0} \hom(\sgm O(n +
m) \otimes_{\SG_n} A^{\otimes n} \otimes M_2 \otimes \dots \otimes
M_m, N)$ such that the canonical diagrams
\[\xymatrix{
  H_A(M_2, \dots, M_m; N) \ar[r]\ar[d] &
  \hom(\sgm O(n + m) \otimes_{\SG_n} A^{\otimes n} \otimes M_2 \otimes \dots
  \otimes M_m, N)
  \ar[d] \\
  {\begin{aligned}
      \hom(\sgm O(k + m) \otimes_{\SG_k} A^{\otimes k}
      \\ \otimes M_2 \otimes \dots \otimes M_m, N)
    \end{aligned}}
  \ar[r] &
  {\begin{aligned}
      \hom((\sgm O(n + m) \otimes \sgm O(k_1) \otimes \dots \sgm
      O(k_{n + m})) \\
      \otimes_{\SG_k} A^{\otimes k} \otimes M_2 \otimes
      \dots \otimes M_m, N)
    \end{aligned}
  }}
\] with $k := \sum_{k = 1}^{m + n} k_i - m$
commute.  There is a natural $A$-operation on the object $H_A(M_2,
\dots, M_m; N)$ which is induced by the canonical morphisms given
componentwise by
\begin{multline*}
  \sgm O(k + 1) \otimes A^{\otimes k} \otimes \hom(\sgm O(n + m) \otimes
  A^{\otimes n} \otimes M_2 \otimes \dots \otimes M_m, N) \\
  \longrightarrow \hom(\sgm
  O(n + m) \otimes A^{\otimes n} \otimes M_2 \otimes \dots \otimes M_m,
  \sgm O(k + 1) \otimes A^{\otimes k} \otimes N) \\
  \longrightarrow \hom(\sgm O(n + m)
  \otimes A^{\otimes n} \otimes M_2 \otimes \dots \otimes M_m, N)
\end{multline*}

\begin{definition}
  The $A$-module $H_A(M_2, \dots, M_m; N)$ as defined above is the
  \emph{lax inner hom from the $A$-modules $M_2, \dots, M_m$ to the
    $A$-module $N$}. 
\end{definition}

\begin{example}
  It is $H_A(; N) = N$ as $A$-modules. 
\end{example}

Let $M_1$ be another $A$-module. The definition of $H_A(M_2, \dots,
M_n; N)$ has been chosen so that the following holds:
\begin{proposition}
  The set of morphisms $P_A(M_1, \dots, M_m) \to N$ of $A$-modules is
  naturally in one-to-one correspondence with the set of morphisms
  $$M_1 \to H_A(M_2, \dots, M_m; N),$$
  i.e.~the functor $H_A(M_2,
  \dots, M_m; \cdot)$ is right adjoint to $P_A(\cdot, M_2, \dots,
  M_m)$. The adjunction morphism is induced by the adjunction morphism
  of the inner hom of $\cat C$.
\end{proposition}
\qed

\begin{example}
  Let $A$ be a unital, associative algebra in $\cat C$, which we also
  view as an $\Ass$-algebra. Let $M_2, \dots, M_n$ and $N$ be modules over
  $A$. Their lax inner hom is given by
  \[
  H_A(M_2, \dots, M_n; N) = \prod_{\sigma \in \SG_n}
  \hom_A(M_{\sigma(1)} \otimes_A \dots \otimes_A M_{\sigma(n)}; N). 
  \]
  where $M_1 := A$ and where the right hand side carries the obvious
  structure of a module over $A$.
\end{example}

\begin{example}
  Let $C$ be a unital, commutative algebra in $\cat C$, which we also
  view as a $\Com$-algebra. Let $M_2, \dots, M_n$ and $N$ be modules
  over $A$. Their lax inner hom is given by
  \[
  H_A(M_2, \dots, M_n; N) = \hom_C(M_2 \otimes_C \otimes \dots
  \otimes_C M_n; N)
  \]
  with the obvious structure of a module over $C$. 
\end{example}

\subsection{Algebras over algebras}

Let $\cat C$ be a cocomplete closed symmetric monoidal category and
let $\sgm O$ be an operad in $\cat C$. Let $A$ be an $\sgm O$-algebra. 

\begin{definition}
  An \emph{$A$-algebra} is a morphism $A \to B$ of $\sgm
  O$-algebras. By abuse of notation, we often write $B$ instead of $A
  \to B$. The \emph{category of $A$-algebras} is the comma category
  $(A, \id)$ where $\id$ is the identity functor on the category of
  $\sgm O$-algebras. 
\end{definition}

\begin{remark}
  Let $A \to B$ be an $A$-algebra. Then $B$ is naturally an $A$-module
  where the $A$-operation is given by the morphisms given
  componentwise by the canonical maps
  \[
  \sgm O(n + 1) \otimes A^{\otimes n} \otimes B \longrightarrow \sgm
  O(n + 1) \otimes B^{\otimes (n + 1)} \longrightarrow B.
  \]
\end{remark}

Let $M$ be an $A$-module. Set 
\[S^*_A(M) := \coprod_{n = 0}^\infty S^n_A(M).\]
There is a unique
equivariant $\sgm O$-multiplication law on $S^*_A(M)$ that is induced by
the morphisms given componentwise by
\begin{multline*}
  \sgm O(n) \otimes (\sgm O(m_1 + k_1) \otimes A^{\otimes k_1} \otimes
  M^{\otimes m_1}) \otimes \dots \otimes (\sgm O(m_n + k_n) \otimes
  A^{\otimes k_n} \otimes M^{\otimes m_n}) \\
  \longrightarrow \sgm O(m + k) \otimes
  A^{\otimes k} \otimes M^{\otimes m}
\end{multline*}
with $k := \sum_{i = 1}^n k_i$ and $m := \sum_{i = 1}^n m_i$. It is
unital, and associative. The canonical morphism $A \to
S^*_A(M)$ given by $A = S^0_A(M)$ is a morphism of $\sgm O$-algebras. 

\begin{definition}
  The \emph{free $A$-algebra over $M$} is the $A$-algebra given by the
  natural morphism $A \to S^*_A(M)$ as defined above. 
\end{definition}

There is a natural morphism of $A$-modules $M \to S^*_A(M)$ induced by
the canonical isomorphism $M = S^1_A(M)$. 

\begin{proposition}
  The free $A$-algebra over $M$ is in fact \emph{free}, i.e.~it is the
  initial object in the comma category $(M, \#)$, where $\#$ is the
  forgetful functor from the category of $A$-algebras to the category
  of $A$-modules. 
\end{proposition}
\qed

\begin{example}
  Let $A$ be a unital, associative algebra, which we also view as an
  $\Ass$-algebra. Let $M$ be an $A$-module. It is
  \[S^*_A(M) = \coprod_{n = 0}^\infty M^{\otimes_A n},\]
  i.e.~the classical tensor algebra.
\end{example}

\begin{example}
  Let $C$ be a unital, commutative algebra, which we also view as a
  $\Com$-algebra. Let $M$ be a $C$-module. It is
  \[S^*_A(M) = \coprod_{n = 0}^\infty M^{\otimes_{A, \SG_n} n},\]
  i.e.~the classical symmetric algebra. 
\end{example}

\begin{example}
  Let $W$ be an object in $\cat C$. The free $A$-algebra
  over $M := F_A(W)$ is given by
  \[S^*_A(M) = U(A) \otimes \coprod_{n = 0}^\infty W^{\otimes_{A, \SG_n}
    n},\] i.e.~the tensor product of the universal envelopping algebra
  and the symmetric algebra over $W$ in $\cat C$.
\end{example}

Let $A \to B$ be an $A$-algebra. Note that every $B$-module is canonically
an $A$-module via the morphism $A \to B$. 

Let $E_1, \dots, E_r$ be $A$-modules. Set
\[S^*_A(M, E_1, \dots, E_r) := \coprod_{n = 0}^\infty S^n_A(M, E_1,
\dots, E_r).\] There is a unique equivariant $S^*_A(M)$-operation on
$S^*_A(M, E_1, \dots, E_r)$ that is induced by the morphisms given
componentwise by
\begin{multline*}
  \sgm O(n + 1) \otimes (\sgm O(m_1 + k_1) \otimes A^{\otimes k_1}
  \otimes M^{\otimes m_1}) \otimes \dots \otimes (\sgm O(m_n + k_n)
  \otimes A^{\otimes k_n} \otimes M^{\otimes m_n}) \\
  \otimes (\sgm O(m_{n
    + 1} + k_{n + 1} + r) \otimes A^{\otimes k_{n + 1}} \otimes M^{m_{n
      + 1}} \otimes E_1 \otimes \dots \otimes E_r)
  \\
  \longrightarrow O(m + k + r) \otimes A^{\otimes k} \otimes
  M^{\otimes m} \otimes E_1 \otimes \dots \otimes E_r
\end{multline*}
with $k := \sum_{i = 1}^{n + 1} k_i$ and $m := \sum_{i = 1}^{n + 1}
m_i$. It is unital, equivariant and associative.

Let $E$ be an $A$-module.
The canonical morphism $E \to S^*_A(M, E)$ induced by
the isomorphism $E = S^1(M, E)$ is a morphism of $A$-modules.

\begin{proposition}
  The $S^*_A(M)$-module $S^*_A(M, E)$ as defined above is a
  \emph{free} object, i.e.~it is the initial object in the comma
  category $(E, \#)$, where $\#$ is the canonical forgetful functor
  from the category of $S^*_A(M)$-modules to the category of
  $A$-modules.
\end{proposition}

\subsection{Connections}

Now we shall propose an extension of the notion of a connection of a
module over a commutative algebra to modules over an algebra of any
operadic type. For this, one has to make use of the lax product as
defined above.

Let $\cat C$ be a cocomplete closed symmetric monoidal additive
category and let $\sgm O$ be an operad in $\cat C$. Let $A$ be an
$\sgm O$-algebra and let $d: A \to M$ be a derivation of $A$. 

Let $E$ be an $A$-module. A morphism $\nabla: E \to P_A(M, E)$ is
\emph{$d$-derivative} if the canonical diagrams
\[\xymatrix@C+96pt{ \sgm O(n + 1) \otimes A^{\otimes n} \otimes E
  \ar[r]^{\substack{\id \otimes (\sum_{p + 1 + q = n} \id^{\otimes p}
      \otimes d \otimes \id^{\otimes q}) \otimes \id \\ + \id \otimes
      \id^{\otimes n} \otimes \nabla}} \ar[d] & {\begin{aligned} \sgm
      O(n + 1) \otimes A^{\otimes (n - 1)} \otimes M \otimes E
      \\
      \oplus \sgm O(n + 1) \otimes A^{\otimes n} \otimes P_A(M, E)
  \end{aligned}}
    \ar[d] \\
  E \ar[r]_{\nabla} & P_A(M, E) }\] commute.

\begin{definition}
  A \emph{$d$-connection on $E$} is a $d$-derivative morphism $E \to
  P_A(M, E)$. 
\end{definition}

\begin{example}
  In this example, we view $A$ as an module over itself.
  The morphism $d: A \to M$ composed with the natural morphism $M \to
  P_A(M, A)$ is a $d$-connection on $A$.
\end{example}

\begin{example}
  Let $A$ be a unital, associative algebra, which we also view as an
  $\Ass$-algebra. Let $d: A \to M$ be an $A$-derivative and $E$ an
  $A$-module. A $d$-connection $\nabla: E \to P_A(M, E) = M \otimes_A
  E \oplus E \otimes_A M$ is a morphism such that the canonical
  diagrams
  \[\xymatrix@C+96pt{
    A \otimes E \ar[r]^{d \otimes \id + \id \otimes \nabla}
    \ar[d] &
    M \otimes E \oplus A \otimes (M \otimes_A E \oplus E \otimes_A M)
    \ar[d] \\
    E \ar[r]_{\nabla} &
    M \otimes_A E \oplus E \otimes_A M
  }\]
  and
  \[\xymatrix@C+96pt{
    E \otimes A \ar[r]^{\id \otimes d + \nabla \otimes \id}
    \ar[d] &
    E \otimes M \oplus (M \otimes_A E \oplus E \otimes_A M) \otimes A
    \ar[d] \\
    E \ar[r]_{\nabla} &
    M \otimes_A E \oplus E \otimes_A M
  }\]
  commute. 
\end{example}

\begin{example}
  Let $C$ be a unital, commutative algebra, which we also view as a
  $\Com$-algebra. Let $d: A \to M$ be a derivation of $A$ and $E$ an
  $A$-module. A $d$-connection $\nabla: E \to P_A(M, E) = M \otimes_A
  E$ is just a $d$-connection in the classical sense.
\end{example}

\begin{example}
  Let $W$ be an object in $\cat C$. There is a canonical
  $d$-connection on $F_A(W) = U(A) \otimes W$ which is induced by the
  morphisms given componentwise by
  \[\xymatrix@C+96pt{
    \sgm O(n + 1) \otimes A^{\otimes n} \otimes W \ar[r]^{\id \otimes
      (\sum_{p + 1 + q = n} \id^{\otimes p} \otimes d \otimes
      \id^{\otimes q}) \otimes \id} & \sgm O(n + 1) \otimes A^{\otimes
      (n - 1)} \otimes M \otimes W \ar[d]
    \\ & P_A(M, F_A(W))}\]
\end{example}

\subsection{Jet modules}

Let $\cat C$ be a cocomplete closed symmetric monoidal abelian
category, and let $\sgm O$ be an operad in $\cat C$. Let $A$ be an
$\sgm O$-algebra and let $d: A \to M$ be a derivation of $A$. 

Let $E$ be an $A$-module and let $\gamma$ be the multiplication map of
$A$ on $E$. Set
\[
J_d E := E \oplus P_A(M, E). 
\]
There is a unique equivariant $A$-operation on $J_d E$, which is
induced by the canonical morphisms given componentwise by
\[\xymatrix@C+78pt{
  \sgm O(n + 1) \otimes A^{\otimes n} \otimes E
  \ar[r]^{\gamma + \id \otimes (\sum_{p + 1 + q = n} \id^{\otimes p}
    \otimes d \otimes \id^{\otimes q}) \otimes \id} &
  {\begin{aligned}
      E \\ \oplus \left(\sgm O(n + 1) \otimes A^{\otimes (n - 1)}
        \otimes M \otimes E\right)
    \end{aligned}}
  \ar[d]
  \\ &
  E \oplus P_A(M, E),
}\]
and
\[\xymatrix@C+24pt{
  \sgm O(n + 1) \otimes A^{\otimes n} \otimes P_A(M, E) \ar[r]^\gamma
  & P_A(M, E).
}\]
This operation is unital, equivariant, and associative. 
\begin{definition}
  The $A$-module $J_d E$ as defined above is the \emph{(one-)$d$-jet
    module of $E$}.
\end{definition}

The $A$-operation on $J_d E$ has been chosen so that the
following proposition holds:
\begin{proposition}
  A morphism $\nabla: E \to P_A(M, E)$ is a $d$-connection if and only
  if $\id + \nabla = (\id, \nabla): E \to E \oplus P_A(M, E)$ is a
  morphism of $A$-modules.
\end{proposition}
\qed

There is a canonical complex
\[\xymatrix{
  0 \ar[r] & P_A(M, E) \ar[r] & J_d E \ar[r] & E \ar[r] & 0, }\] of
$A$-modules where the second morphism is the canonical inclusion and
the third morphism is the canonical projection. This sequence is exact
as the underlying sequence of objects in $\cat C$ is exact.
\begin{definition}
  The short exact sequence
  \[\xymatrix{
    0 \ar[r] &
    P_A(M, E) \ar[r] &
    J_d E \ar[r] &
    E \ar[r] &
    0
  }\]
  of $A$-modules is the \emph{(first) $d$-jet module sequence of
    $E$}. 
\end{definition}
For the classical (i.e.~non-operadic case) in the geometrical setting
this sequence is discussed in~\cite{Kapranov99}.

\section{Model categories}

The rest of the article uses the language of model categories to do
the necessary homological algebra. In order to fix the notation, we
repeat the main notions.

\subsection{Retracts}

Let $\cat C$ be a category. Let $f$ and $g$ be two morphisms in $\cat C$ 
such that a commutative diagram of the form
\[\xymatrix{
  A \ar@{=}@/^/@<1ex>[rr] \ar[d]_f\ar[r] & C \ar[d]_g\ar[r] & A \ar[d]^f \\
  B \ar@{=}@/_/@<-1ex>[rr] \ar[r] & D \ar[r] & B
}\]
exists.
\begin{definition}
  In this situation, $f$ is a \emph{retract of $g$}. 
\end{definition}

Let $\cat W$ a subcategory of $\cat C$. 
\begin{definition}
  The \emph{subcategory $\cat W$ is closed under retracts} if for any
  two morphisms $f$ and $g$ such that $f$ is a retract of $g$ and $g$
  is in $\cat W$, $f$ is as well.
\end{definition}

\subsection{Lifting properties}
Let $\cat C$ be a category. Consider commutative diagrams of the form
\[\xymatrix{
  A \ar[r]^f \ar[d]_i & X \ar[d]^p \\
  B \ar[r]_g \ar@{-->}[ur]^h & Y
}\]
in $\cat C$. 
\begin{definition}
  A morphism $i$ in $\cat C$ has the \emph{left lifting property with
    respect to a morphism $p$ in $\cat C$} (and $p$ the \emph{right
    lifting property with repect to $i$}) if for all morphism $f$ and
  $g$ in $\cat C$ such that the solid square in the above diagram
  commutes, there exists a morphism $h$ making the whole diagram
  commutative. 
\end{definition}

\subsection{Two-out-of-three axiom}

Let $\cat C$ be a category and $\cat W$ a subcategory of $\cat C$. 

\begin{definition}
  The \emph{subcategory $\cat W$ satisfies the two-out-of-three axiom}
  if for any three morphisms $f$, $g$ and $h$ of $\cat C$ with $h = g
  f$ such that two of them are in $\cat W$, the third is as well.
\end{definition}

\subsection{Model categories}

Let $\cat C$ be a category. Let us assume that there are three
distinguished subcategories of $\cat C$ whose morphisms are called
\emph{weak equivalences}, \emph{cofibrations} and \emph{fibrations},
respectively.

Morphisms that are weak equivalences and cofibrations at the same time
are \emph{acyclic (or trivial) cofibrations}. Morphisms that are weak
equivalences and fibrations at the same time are \emph{acyclic (or
  trivial) fibrations}.

\begin{definition}
  The three distinguished subcategories define a \emph{model structure
    on $\cat C$} if 
  \begin{itemize}
  \item
    the subcategory of the weak equivalences
    fulfills the two-out-of-three axiom,
  \item each of the three subcategories of weak equivalences,
    cofibrations, and fibrations is closed under retracts,
  \item
    the trivial cofibrations have the left lifting property
    with respect to the fibrations, the trivial fibrations have the
    right lifting property with respect to the cofibration, and
  \item each
    morphisms $f$ in $\cat C$ has two functorial factorisations $f = p
    \circ i = q \circ j$ such that $p$ is a fibration, $i$ is a acyclic
    cofibration, $q$ is a acyclic fibration and $j$ is a cofibration.
  \end{itemize}
\end{definition}

\begin{definition}
  A bicomplete category $\cat C$ together with a model structure on it
  is a \emph{model category}.
\end{definition}

\begin{definition}
  An object in a model category is \emph{cofibrant} if the unique
  morphism from the initial object to the object is a cofibration. An
  object is \emph{fibrant} if the unique morphism from the object to
  the terminal object is a fibration.
\end{definition}

For more on model categories, we refer the reader to the
monograph~\cite{Hovey99} and the references therein.

\subsection{Monoidal model categories}

Let $\cat C$ be a bicomplete closed (symmetric) monoidal category. 

Let $f: A \to B$ and $g: X \to Y$ be two morphisms in $\cat C$. Let
\[f \boxempty g: A \otimes Y \amalg_{A \otimes X} B \otimes X \to B \otimes
Y\]
be the natural morphism induced by $f$ and $g$. 
\begin{definition}
  The morphism $f \boxempty g$ as defined above is the \emph{pushout tensor
    product of $f$ and $g$}. 
\end{definition}

Assume that $\cat C$ is endowed with a model structure. Let
\[\xymatrix{0 \ar[r] & Q\unit \ar[r]^q & \unit}\]
be the functorial factorisation of $0 \to \unit$ into a cofibration
followed by an acyclic fibration. 
\begin{definition}
  The category $\cat C$ is a \emph{symmetric monoidal model category}
  if for any two cofibrations $f$ and $g$, their pushout tensor
  product $f \boxempty g$ is a cofibration which is acyclic
  if $f$ or $g$ is acyclic, and if the natural morphisms $q \boxempty
  X: Q\unit \otimes X \to \unit \otimes X$ are weak equivalences for
  all cofibrant objects $X$. 
\end{definition}
See also~\cite{Hovey98} and~\cite{SchwedeShipley00}.

\section{Cofibrantly generated model categories}

Generally it is not an easy task to construct model structures on
categories ``by hand''. There is, however, a general method to
construct certain model structures, the so-called cofibrantly
generated model structures. This method is based on Quillen's ``small
object argument'', see~\cite{Quillen67} and~\cite{Hovey99}. We repeat
the main notions here.

\subsection{Smallness}

A \emph{limit ordinal} is an ordinal that is not the direct successor
of an ordinal. 

\begin{definition}
  Let $\lambda$ be a limit ordinal. The \emph{cofinality $\cofin
    \lambda$ of $\lambda$} is the least cardinal $\kappa$ such that
  there exists a subset $T$ of $\lambda$ with $\card T = \kappa$ and
  $\sup T = \lambda$.
\end{definition}

\begin{example}
  It is $\cofin \kappa = \kappa$ for each cardinal $\kappa$.
\end{example}

Let $\cat C$ be a cocomplete category, $A$ an object in $\cat C$ and
$\kappa$ a cardinal. Let $\cat D$ be a subcategory of $\cat C$.
\begin{definition}
  The object $A$ is \emph{$\kappa$-small (relative to $\cat D$)} if,
  for every ordinal $\lambda$ with $\cofin \lambda > \kappa$ and every
  colimit-preserving functor $X: \lambda \to \cat C$ ($X: \lambda \to
  \cat D$), the natural map $\colim_{\mu < \lambda} \hom(A, X_\mu) \to
  \hom(A, \colim_{\mu < \lambda} X_\mu)$ is an isomorphism.
\end{definition}

\begin{example}
  If $A$ is $\kappa$-small for a cardinal $\kappa$, it is also
  $\kappa'$-small for each cardinal $\kappa' \ge \kappa$. 
\end{example}

\begin{definition}
  The object $A$ is \emph{small} if it is $\kappa$-small for some
  cardinal $\kappa$. 
\end{definition}

\begin{theorem}
  Let $\cat C$ be a Grothendieck category (i.e.~a bicomplete abelian
  category with a generator and exact filtered colimits). Any object
  of $\cat C$ is small, i.e.~small relative to the whole category
  $\cat C$. 
\end{theorem}

\begin{proof}
  Usually, the proof is given by using an embedding theorem for
  Grothendieck categories or by using that any Grothendieck category
  is locally presentable. For one these proofs, we refer to~\cite{Hovey01}.
  
  In the appendix, we give a proof that does not use any deep theorem
  and relies solely on the basic properties of a Grothendieck
  category.
\end{proof}

\subsection{Cells}

Let $\cat C$ be a cocomplete category. Let $I$ be a class of morphisms
in $\cat C$. 

A morphism in $\cat C$ is \emph{$I$-injective} if it has the right
lifting property with respect to every morphism in $I$. 

\begin{definition}
  The class of all $I$-injective morphisms is denoted by $I$-inj. 
\end{definition}

For example, an object $Y$ in $\cat C$ is injective (in the sense of
injective modules) if $0 \to Y$ is in $I$-inj, where $I$ is the class
of monomorphisms in $\cat C$.

A morphism in $\cat C$ is an \emph{$I$-cofibration} if it has the left
lifting property with respect to every morphism in $I$-inj. 

\begin{definition}
  The class of all $I$-cofibrations is denoted by $I$-cof. 
\end{definition}

In particular, it is $I$ a subset of $I$-cof.

A morphism in $\cat C$ is a \emph{relative $I$-cell complex} if it is
a transfinite composition of pushouts of morphisms in $I$. Here, a
transfinite composition is just a colimit over a colimit-preserving
functor $X: \lambda \to \cat C$, where $\lambda$ is a limit
ordinal. The monograph~\cite{Hovey99} is also a good reference for
this notion.

\begin{definition}
  The class of all relative $I$-cell complexes is denoted by $I$-cell. 
\end{definition}

For example, relative CW-complexes in topology are exactly the
relative $I$-cell complexes when $I$ is the set of all morphisms
$\partial \set D^n \to \set D^n$, where $\set D^n$ denotes the
\emph{$n$-disk}.

\begin{remark}
  It is $I$-cell a subclass of $I$-cof.
\end{remark}

\subsection{Cofibrantly generated model categories}

Our main reference for this subsection is again~\cite{Hovey99}.

Let $\cat C$ be a model category. Let $I$ and $J$ be two sets of
morphisms in $\cat C$. 
\begin{definition}
  The model category $\cat C$ is \emph{cofibrantly generated} if the
  domains of the morphisms in $I$ are small relative to $I$-cell, the
  domains of the morphisms in $J$ are small relative to $J$-cell, the
  class of fibrations is $J$-inj, and the class of trivial fibrations
  is $I$-inj. $I$ is the \emph{set of generating cofibrations} and $J$
  the \emph{set of generating acyclic cofibrations}.
\end{definition}

Let $R$ be a unital commutative ring (in the ordinary sense). 
\begin{example}
  The category of cochain complexes of $R$-modules is a cofibrantly
  generated model category whose weak equivalences are the
  quasiisomorphisms and whose fibrations are all degree-wise
  surjective morphisms. For the proof and the generating (acyclic)
  cofibrations see, e.g.,~\cite{Hovey99}.

  In fact, the category is a symmetric monoidal model category. 
\end{example}

Let $(X, \mathcal O_X)$ be a ringed space with finite hereditary
global dimension (see~\cite{Hovey01} in this context), e.g.~$X$ is a
finite-dimensional noetherian topological space or a
finite-dimensional locally compact Hausdorff space that is countably
at infinity.
\begin{example}
  This example is one of the main results in~\cite{Hovey01}.  The
  category of cochain complexes of $\mathcal O_X$-modules is a
  cofibrantly generated model category whose weak equivalences are the
  quasiisomorphisms and whose fibrations are the degree-wise
  surjections with degree-wise flabby kernel.

  In fact, the category is a symmetric monoidal model category. This
  example subsumes the previous one (for $X = \{*\}$, the one-pointed
  space).
\end{example}

\subsection{Categories of modules}

This subsection builts on the ideas and results
in~\cite{SchwedeShipley00} and~\cite{Hovey98}.

Let $\cat C$ be a cofibrantly generated symmetric monoidal model
category. Let $U$ be a monoid (i.e.~a unital, associative algebra) in
$\cat C$.

Let $I$ be the set of generating cofibrations and $J$ be the set of
generating trivial cofibrations. 

Assume that the following conditions hold: The domains of the
morphisms in $I$ are small relative to $(U \otimes I)$-cell. The
domains of the morphisms in $J$ are small relative to $(U \otimes
J)$-cell. Every morphism in $(U \otimes J)$-cell is a weak
equivalence. (If these conditions hold, we say the cofibrantly
generated monoidal model category $\cat C$ \emph{fulfills the monoid
  axiom for $U$}.)

The following theorem is proved in~\cite{Hovey98}.
\begin{theorem}
  There is a cofibrantly generated model structure on the category of
  left $U$-modules, where a morphism of left $U$-modules is a weak
  equivalence (respectively fibration) if and only if it is a weak
  equivalence (respectively fibration) in $\cat C$. The set of
  generating cofibrations is $U \otimes I$, the set of
  generating acyclic cofibrations is $U \otimes J$.
\end{theorem}
\qed

We should remark that Hinich provides us in~\cite{Hinich97} with a
different proof for the existence of a model structure on a category
of modules over an algebra.

\begin{definition}
  In the situation of the previous theorem, we say that the category
  of left $U$-modules \emph{admits a model structure over $\cat C$}.
\end{definition}

Let $\cat C$ be a cofibrantly generated symmetric monoidal abelian
model category that is a Grothendieck category. Assume that $\cat C$
fulfills the monoid axiom. Let $\mathcal O$ be an operad in $\cat C$
and let $A$ be an algebra over $\mathcal O$.
\begin{example}
  The category of modules over $A$ admits a model structure over $\cat
  C$ as it is isomorphic to the category of left $U(A)$-modules. 
\end{example}

\subsection{The category of modules over an algebra}

In this subsection, we apply the previous theorem to the case of
modules over an algebra over an operad.

Let $\cat C$ be a cofibrantly generated closed symmetric model
category. Let $I$ be the set of generating cofibrations and $J$ be the
set of generating acyclic cofibrations of $\cat C$. Let $\sgm O$ be an
operad in $\cat C$. Let $A$ be an $\sgm O$-algebra.

Consider the monoid $U(A)$ in $\cat C$. Assume that $\cat C$ satisfies
the monoid axiom for $U(A)$. Recall that the free module functor $F_A$
is just tensoring by $U(A)$.
\begin{example}
  The category of $A$-modules becomes cofibrantly generated model
  category, where a morphism of $A$-modules is a weak equivalence
  (respectively a fibration) if and only if it is a weak equivalence
  (respectively a fibration) in $\cat C$. The set of generating
  cofibrations is $F_A(I)$, the set of generating acyclic cofibrations
  is $F_A(J)$.
\end{example}

\subsection{The category of modules over a ringed space as an example}

Let $(X, \mathcal O_X)$ be a ringed space with finite hereditary
global dimension. The tensor product over $\mathcal O_X$ makes the
category of cochain complexes of $\mathcal O_X$-modules a closed
symmetric monoidal category, which is by a result of Hovey
(see~\cite{Hovey01}) in fact a closed symmetric monoidal model
category.

\begin{proposition}
  The category of $\mathcal O_X$-modules satisfies the monoid axiom
  for any monoid.
\end{proposition}

\begin{proof}
  Let $J$ be the set of generating acyclic cofibrations of the
  category $\cat X^*$ of cochain complexes of $\mathcal O_X$-modules.
  We have to show that transfinite compositions of pushouts of
  morphisms in $U \otimes J$ are quasiisomorphism for any object $U$
  in $\cat X^*$. As cohomology commutes with filtered colimits, we
  just have to consider pushouts of morphisms in $U \otimes J$. That
  these are quasiisomorphisms can be verified stalk-wise. However,
  stalk-wise, the morphisms in $J$ (which are given
  in~\cite{Hovey01}) are by definition either isomorphisms of exact
  complexes of free modules or inclusions of the zero complex into
  exact complexes of free modules.  Thus stalk-wise, the morphisms in
  $U \otimes J$ are either isomorphisms or inclusions of the zero
  complex into exact complexes. Pushouts of such morphisms are
  quasiisomorphisms.
\end{proof}   

Let $U$ be any monoid in the category of $\mathcal
O_X$-modules. As any object in the Grothendieck category of cochain
complexes of $\mathcal O_X$-modules is small and this category
satisfies the monoid axiom, we arrive at the following example due
to the theorem in the previous subsection:
\begin{example}
  The category of (left) $U$-modules is a cofibrantly generated model
  category, where a morphism of left $U$-modules is a weak
  equivalence if and only if it is a quasiisomorphism, and where a
  morphism of left $U$-modules is a fibration if and only if it is a
  degree-wise surjection with degree-wise flabby kernel. 
\end{example}

\subsection{Fibrant replacements}

Let $\cat C$ be a model category with terminal object $*$.

By the axioms of a model structure, every morphism $f: X \to *$ in
$\cat C$ can be decomposed functorially as $f = p i$ where $i: X \to
RX$ is an acyclic cofibration and $p: RX \to *$ is a fibration. We
call $RX$ a \emph{(the) fibrant replacement for $X$}. This leads to
the following definition:

\begin{definition}
  A \emph{fibrant replacement functor $R: \cat C \to \cat C$} is part
  of a natural transformation $\eta: \id \Rightarrow R$ such that each
  natural morphism $X \to RX$ is a weak equivalence and each $RX$ is a
  fibrant object for all objects $X$ in $\cat C$ 
\end{definition}

\begin{example}
  As we have seen above, each model category comes with a canonical
  fibrant replacement functor. 
\end{example}

The main reason why one considers other fibrant replacement functors
is the presence of a tensor product and the question of
compatibility
 
Assume that $\cat C$ is a closed symmetric monoidal
model category with a fibrant replacement functor $R$. 

The following definition is from~\cite{BergerMoerdijk03}.
\begin{definition}
  The fibrant replacement functor $R$ is \emph{symmetric monoidal} if
  it is symmetric monoidal as a functor of symmetric monoidal
  categories and the natural diagrams
  \[\xymatrix{
    X \otimes Y \ar[d] \ar[dr] \\
    RX \otimes RY \ar[r] & R(X \otimes Y),
  }\]
  natural in $X$ and $Y$, commute. 
\end{definition}

\begin{remark}
  By adjunction, every symmetric monoidal fibrant replacement functor
  $R$ in $\cat C$ defines morphisms $\hom(X, X) \to \hom(RX, RX)$,
  natural in $X$.
\end{remark}

The following proposition yields plenty of examples of closed
symmetric monoidal categories with a fibrant replacement functor $R$. 
\begin{proposition}
  Let $\cat C$ be cofibrantly generated with $J$ being the set of
  generating acyclic cofibrations. Assume that $J \boxempty J$ is a subset
  of $J$-cell. Then $\cat C$ admits a symmetric monoidal fibrant
  replacement functor.
\end{proposition}

\begin{proof}
  The proof makes use of Quillen's small object argument. In fact, we
  start with the standard proof of constructing a functorial
  factorisation into a acyclic cofibration and a fibration.
  
  Let $J$ be the set of generating acyclic cofibrations. Let $\kappa$
  be a cardinal such that every domain of the morphisms in $J$ is
  $\kappa$-small. Let $\lambda$ be a limit ordinal with $\cofin
  \lambda > \kappa$. 
  
  For each object $X$ in $\cat C$, we shall naturally define a functor
  $(X_\beta)_{\beta < \lambda}: \lambda \to \cat C$ inductively as
  follows: We set $X_0 := X$. For each ordinal $\beta < \lambda$ let
  $X_{\beta + 1}$ be the simultaneous pushout of all diagrams of the
  form
  \[\xymatrix{
    A \ar[d] \ar[r] & X_\beta \ar@{-->}[d] \\
    B \ar@{-->}[r] & X_{\beta + 1},
  }\]
  where the vertical arrow runs through all morphisms in $J$ and the
  horizontal morphism through all morphisms into $X_\beta$. For a
  limit ordinal $\beta < \lambda$, we set $X_{\beta} :=
  \colim_{\beta' < \beta} X_{\beta}$. Finally set $RX := \colim_{\beta
    < \lambda} X_\beta$. 
  
  The natural morphism $X \to RX$ is in $J$-cell and thus a weak
  equivalence by the axioms of a cofibrantly generated model category.
  We have to show that the morphism $RX \to *$ has the right lifting
  property with respect to morphism in $J$-cell, which means that $RX
  \to *$ is a fibration. (In fact, we shall construct natural lifts.)
  To show this, consider a pushout diagram
  \[\xymatrix{
    A \ar[r]\ar[d]_j & Z \ar[d] \\
    B \ar[r] & Z',
  }\]
  with $j \in J$. Assume that there is given a morphism $Z \to RX$. We
  want to extend this map to $Z'$. The composition $A \to Z \to RX$
  factors by the $\kappa$-smallness of $A$ through some $X_\beta$.  By
  construction, this morphism $A \to X_\beta$ extends naturally to a
  morphism $B \to X_{\beta + 1}$. Passing to the colimit shows that
  the morphism $A \to RX$ extends naturally to a morphism $B \to RX$. By
  definition of the pushout, this defines an extension $Z' \to RX$ of
  $Z \to RX$. Now, every morphism in $J$-cell is a filtered colimit of
  such pushouts $Z \to Z'$. Thus $RX \to *$ has the right lifting
  property with repect to morphisms in $J$-cell. 
  
  It remains to construct natural morphisms $RX \otimes RY \to R(X
  \otimes Y)$. By construction and the fact that pushout products of
  maps in $J$ are in $J$-cell, the morphism $X \otimes Y \to RX
  \otimes RY$ is in $J$-cell, so by the above considerations, there
  exists a natural lift $RX \otimes RY \to R(X \otimes Y)$ making the
  diagram
  \[\xymatrix{
    X \otimes Y \ar[r] \ar[d] & R(X \otimes Y) \\
    RX \otimes RY \ar@{-->}[ur]
  }\]
  commutative.
\end{proof}

\begin{remark}
  Assume that pushouts and filtered colimits of monomorphisms are
  monomorphisms. Then the symmetric monoidal fibrant replacement
  functor constructed in the proof above maps monomorphisms to
  monomorphisms. 
\end{remark}

For example, using this proposition, the symmetric monoidal model
category of cochain complexes of $\mathcal O_X$-modules over a
topological space of finite hereditary global dimension admits a
symmetric monoidal fibrant replacement functor.

The remark holds for example true in the category of cochain complexes
$\mathcal O_X$-modules considered in the previous paragraph.

\subsection{The homotopy category}

Let $\cat C$ be a model category. By $\Ho \cat C$ we denote the
``category'' which is the localisation of $\cat C$ by all weak
equivalences, i.e.~we formally invert all weak equivalences in $\cat
C$. The word ``category'' is printed in quotes as it is a priori not
clear if $\Ho \cat C$ is really a (locally small) category, i.e.~if
each class of morphisms between two objects is in fact a set. In this
section, we give one argument for the well-known fact that this is
true.

Let $s: X' \to X$ be a weak equivalence and $f: X' \to Y$ be any
morphism in $\cat C$. This defines a morphism $f s^{-1}: X \to Y$ in
$\Ho \cat C$. 

Let $QX \to X$ be a cofibrant replacement for $X$, i.e.~$QX$ is a
cofibrant object and $QX \to X$ is an acyclic fibration. Let $Y \to
RY$ be a fibrant replacement for $Y$, i.e.~$RY$ is a fibrant object
and $Y \to RY$ is an acyclic cofibration. Further let us decompose the
morphism $s: X' \to X$ as $s = p \circ i$ where $i: X' \to \tilde X$
is an acyclic cofibration and $p: \tilde X \to X$ is an acyclic
fibration. By functoriality of the fibrant replacement functor, the
morphism $f: X' \to Y$ induces a morphism $Rf: RX' \to RY$ where $RX'$
is the fibrant replacement of $X'$.  By two lifting properties due to
the axioms of a model structure there exist dashed arrows making the
diagram
\[\xymatrix{
  & & RX' \ar[r]^{Rf} & RY \\
  & \tilde X \ar@{-->}[ur] \ar[d] & X' \ar[u] \ar[l] \ar[dl]^s
  \ar[r]_f & Y \ar[u] \\
  QX \ar@{-->}[ur] \ar[r] & X }\] commutative. By this construction,
the morphism $f s^{-1}$ can be written as the composition of the
formal inverse of $Y \to RY$, a morphism $QX \to RY$ in $\cat C$ (!),
and the formal inverse of $QX \to X$. We call the result of this
construction a \emph{standard representation of the morphism $f s^{-1}$ in
  $\Ho \cat C$}. Generally, we shall call a \emph{standard
  representation of $f s^{-1}$} a representation of $f s^{-1}$ of the
form
\[X \longrightarrow Q^M X \longrightarrow R^N Y \longrightarrow Y,\]
where the first arrow and third arrow are given by the inverses (in
$\Ho \cat C$!) of the arrows given by the cofibrant and fibrant
replacement and the second arrow is a morphism in $\cat C$.

\begin{lemma}
  Let $t: Y \to Y'$ be a weak equivalence and $g: X \to Y'$ any
  morphism in $\cat C$. Assume that $Y$ is fibrant. Then there exists
  a weak equivalence $s: \tilde X \to X$ and a morphism $f: \tilde X
  \to Y$ such that $t^{-1} g = f s^{-1}$ in $\Ho \cat C$.
\end{lemma}

\begin{proof}
  Consider the pullback diagram
  
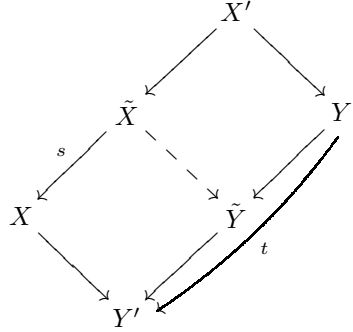
\[\xymatrix{
    & & X' \ar[dl] \ar[dr] \\
    & \tilde X \ar[dl]_s \ar@{-->}[dr] & & Y \ar[dl] \ar@/^/@<1ex>[ddll]^t \\
    X \ar[dr] & & \tilde Y \ar[dl] \\
    & Y' }\] made up from the outer solid arrows, where $X' \to \tilde
  X$ is a cofibration, $Y \to \tilde Y$ is an acyclic cofibration,
  $\tilde X \to X$ and $\tilde Y \to Y'$ are acyclic fibrations and
  where the composition $Y \to Y'$ is the given morphism $t$. The
  dashed line making the diagram commutative exists due to the left
  lifting property of cofibrations with respect to acyclic fibrations.
  Let $s$ be the acyclic fibration given by $\tilde X \to X$. As $Y$
  is fibrant by assumption, the left lifting property of acyclic
  cofibrations with respect to fibrations yields an left inverse $j: Y
  \to \tilde Y$ of the acyclic cofibration $\tilde Y \to Y$.  Finally
  let $f: \tilde X \to Y$ be the composition of $\tilde X \to \tilde
  Y$ with this inverse $j$.
\end{proof}

\begin{lemma}
  Let $QX \to X$ be a cofibrant replacement of an object $X$ and $Y
  \to RY$ be a fibrant replacement of an object $Y$ in $\cat C$. Let
  $X \to Y$ be a morphism in $\Ho \cat C$. Then it can be written as
  a composition of the formal inverse of $QX \to X$, a proper morphism
  $QX \to RY$ in $\cat C$, and the formal inverse of $Y \to RY$. 
\end{lemma}

\begin{proof}
  Any morphism $QX \to RY$ in $\Ho \cat C$ can be written as a
  composition $f_n \circ s_n^{-1} \circ \dots \circ f_1 \circ
  s_1^{-1}$ where the $f_i$ are proper morphisms in $\cat C$ and the
  $s_i$ are weak equivalences in $\cat C$. We may assume that all
  intermediate objects are fibrant. By the previous lemma, this
  morphism can also be written as a composition $f \circ s^{-1}$ where
  $f$ is any morphism in $\cat C$ and $s$ is a weak equivalence in
  $\cat C$. Finally by the construction above, the composition $f
  \circ s^{-1}$ is given by a single morphism $\tilde f: QX \to RY$ in
  $\cat C$.
\end{proof}

Thus $\Ho \cat C$ is in fact a category. 
\begin{definition}
  The category $\Ho \cat C$ is the \emph{homotopy category of $\cat
    C$}. 
\end{definition}

\section{The operadic Atiyah class}

\subsection{Cofibrations}

Let $\cat C$ be a cofibrantly generated closed symmetric monoidal
model category that satisfied the monoid axiom for any monoid. Let
$\sgm O$ be an operad in $\cat C$. Let $A$ be an $\sgm O$-algebra.

Let $f_1: M_1 \to N_1, \dots, f_n: M_n \to N_n$ be morphisms of
$A$-modules. Let
\begin{multline*}
  P_{A, \boxempty}(f_1, \dots, f_n): P_A(N_1, M_2, \dots, M_n)
  \amalg_{P_A(M_1,
    \dots, M_n)} P_A(M_1, N_2, \dots, N_n) \\
  \longrightarrow P_A(N_1, \dots, N_n)
\end{multline*}
be the induced morphism of $A$-modules. 

Let $g: X \to Y$ be another morphism of $A$-modules. Let
\begin{multline*}
  H_{A, \boxempty}(f_2, \dots, f_n, g): H_A(N_2, \dots, N_n; R) \\
  \longrightarrow
  H_A(M_2, \dots, M_n; R) \times_{H(M_2, \dots, M_n; S)} H_A(N_2, \dots,
  N_n; S)
\end{multline*}
be the canonical morphism of $A$-modules. By adjunction there is a
one-to-one correspondence between diagrams of $A$-modules of the form
\[\xymatrix{
  P_A(M_1, N_2, \dots, N_n) \amalg_{P_A(M_1, \dots, M_n)} P_A(N_1, M_2,
  \dots, M_n) \ar[r] \ar[d] &
  R \ar[d] \\
  P_A(N_1, \dots, N_n) \ar@{-->}[ur] \ar[r] & S
}\]
and
\[\xymatrix{
  M_1 \ar[r]\ar[d] & H_A(N_2, \dots, N_m; R) \ar[d] \\
  N_1 \ar[r]\ar@{-->}[ur] &
  H_A(M_2, \dots, M_n; R) \times_{H_A(M_2, \dots, M_n; S)} H_A(N_2, \dots,
  N_n; S). 
}\]

The following lemma has been inspired by~\cite{Hovey01}. The
proof is based on an idea from that article.
\begin{lemma}
  Let $I_1, \dots, I_n$ and $K$ classes of morphisms of
  $A$-modules. Assume that for all $i_i \in I_i$, it is
  $P_{A, \boxempty}(i_1, \dots, i_n)$ a morphism in $K$. Then
  $P_{A, \boxempty}(f_1, \dots, f_n)$ is a morphism in $K$-cof whenever
  each $f_i$ is a morphism in $I_i$-cof. 
\end{lemma}

\begin{proof}
  By assumption and adjointness, the morphisms in $I_1$ have the left
  lifting property with repects to morphisms of the form
  $H_{A, \boxempty}(i_2, \dots, i_n; k)$ with $i_j \in I_j$ where $k$ is
  in $K$-inj. Thus every morphism $f_1$ in $I_1$-cof has the left
  lifting property with repects to morphisms of the form
  $H_{A, \boxempty}(i_2, \dots, i_n; k)$ where $k$ is in $K$-inj.
  Applying adjointness again, we see that morphisms of the form
  $P_{A, \boxempty}(f_1, i_2, \dots, i_n)$ are morphisms in $K$-cof. Then
  one repeats analogous arguments to replace successively each
  $i_j$ by an $f_j$ in $I_j$-cof, $j \ge 2$.
\end{proof}

Recall the free module functor $F_A$.
\begin{lemma}
  Let $I$ be the set of generating cofibrations of $\cat C$ and $J$ be
  the set generating acyclc cofibrations of $\cat C$. Without loss of
  generality, we may assume that $J \subset I$.  Then
  $P_{A, \boxempty}(F_A(i_1), \dots, F_A(i_n))$ is a cofibration of
  $A$-modules for all $i_j \in I$ which is acyclic if there is a $k
  \in \{1, \dots, n\}$ with $i_k \in J$.
\end{lemma}

\begin{proof}
  It is $P_{A, \boxempty}(F_A(i_1), \dots, F_A(i_n)) = F_A(i_1
  \boxempty (i_2 \otimes \dots \otimes i_n))$. Now use that $\cat C$
  is a monoidal model category that satisfies the monoid axiom. 
\end{proof}

Recall the generating (acyclic) cofibrations of the model structure on
the category of $A$-modules.
\begin{proposition}
  Let $f_1, \dots, f_n$ be cofibrations of $A$-modules. Then
  $P_{A, \boxempty(f_1, \dots, f_n)}$ is a cofibration which is acyclic
  if one of the $f_i$ is acyclic. 
\end{proposition}

\begin{proof}
  This follows directly from putting the previous two lemmata
  together.
\end{proof}

\subsection{Fibrations}

Let $\cat C$ be a cofibrantly generated closed symmetric model
category that satisfies the monoid axiom and has a symmetric monoidal
fibrant replacement functor $R: \cat C \to \cat C$. Let $\sgm O$ be an
operad in $\cat C$. Let $A$ be an $\sgm O$-algebra. 

Let $M$ be a module over $A$. The $A$-operation on $M$ is given by a
morphism $U(A) \to \hom(M, M)$. By composition with the natural
morphism $\hom(M, M) \to \hom(RM, RM)$, this defines an $A$-operation
on $RM$. This makes $RM$ naturally an $A$-module, i.e.~the functor $R$
on $\cat C$ induces an endofunctor $R$ on the category of $A$-modules,
also equipped with a natural transformation from the identity
endofunctor to itself.

\begin{proposition}
  The so-defined functor $R$ from the category of $A$-modules into
  itself is a fibrant replacement functor on the model category of
  $A$-modules. 
\end{proposition}

\begin{proof}
  This follows at once from the fact that the forgetful functor from
  the category of $A$-modules to $\cat C$ reflects fibrations and weak
  equivalences.
\end{proof}

In a certain sense, this fibrant replacement functor is compatible
with the lax tensor products: Let $M_1, \dots, M_n$ be $A$-modules. 

\begin{remark}
  As the fibrant replacement functor $R: \cat C \to \cat C$ is
  symmetric monoidal, there exist natural commutative diagrams of the
  form
  \[\xymatrix{
    P_A(M_1, \dots, M_n) \ar[dr]\ar[d] \\
    P_A(R M_1, \dots, RM_n) \ar[r] & R P_A(M_1, \dots, M_n)
  }\]
  of $A$-modules. Moreover, these morphisms are compatible with the
  coherence diagrams involving the functors $P_A$ that have been
  mentioned in the previous section. 
\end{remark}

\subsection{Symmetric monoidal model categories of cochain complexes}

Let $\cat C$ be a bicomplete closed symmetric monoidal Grothendieck
category. Let the category of cochain complexes $\cat C^*$ over $\cat
C$ be endowed with the structure of a monoidal model category. We make
the following further assumptions:
\begin{itemize}
\item The weak equivalences in $\cat C^*$ are exactly the
  quasiisomorphisms. 
\item The monoidal model structure on $\cat C^*$ satisfies the monoid
  axiom for any monoid.
\item The property of a morphism between two cochain complexes being a
  fibration can be tested degree-wise. 
\item There exists a symmetric monoidal fibrant replacement functor
  $R: \cat C^* \to \cat C^*$. 
\item The translation functor on the category of cochain complexes
  $\cat C^*$ respects all of the structure. 
\end{itemize}

\begin{definition}
  Under these assumptions, we call the category $\cat C^*$ is a
  \emph{symmetric monoidal model category of cochain complexes} for
  short.
\end{definition}

As our own interest lies in examples coming from geometry, the
following is important for us:
\begin{example}
  Let $(X, \mathcal O_X)$ be a ringed space with finite hereditary
  dimension. The category cochain complexes of $\mathcal
  O_X$-modules is canonically a symmetric monoidal category of cochain
  complexes. 
\end{example}

\begin{remark}
  The axioms imposed on the symmetric monoidal model category of
  cochain complexes on $\cat C^*$ imply in particular that there is a
  symmetric monoidal model structure on the category of modules over
  an algebra over an operad in $\cat C^*$ that possesses a fibrant
  replacement functor that is compatible with the lax tensor
  products. 
\end{remark}

\subsection{Extension classes}

Let $\cat C$ be a symmetric monoidal model category of cochain complexes. 

Let
\[\xymatrix{
  0 \ar[r] & E' \ar[r]^f & E \ar[r] & E'' \ar[r] & 0 }\] be a short
exact sequence of objects in $\cat C^*$. Recall the definition of the
\emph{cone $\cone f$ of $f$}: It is $(\cone f)^n = E^{\prime n + 1}
\oplus E^{\prime n}$ and the differential on $(\cone f)^n$ is given by
$- d_{E'}^{n + 1} + f + d_E^n$. The morphism $E \to E''$ induces (via
the projection $\cone f \to E$, which is itself not a map of
complexes) a quasiisomorphism $s: \cone f \to E''$. There is further a
canonical projection morphism $p: \cone f \to E'[1]$.

Let $R: \cat C \to \cat C$ be a fibrant replacement functor. Then
$\cone (Rf)$ is a fibrant object. 

\begin{definition}
  The \emph{extension class associated to the short exact sequence
    above} is the morphism $p s^{-1}: E'' \to E'[1]$ in $\Ho \cat C$.
\end{definition}

Recall the standard representation of the morphism $p s^{-1}$: Let
$QE'' \to E''$ be a cofibrant replacement of $E''$ and $E' \to RE'$ be
a fibrant replacement of $E'$ (which makes $E'[1] \to RE'[1]$ a
fibrant replacement of $E'[1]$). Then dashed arrows exist making the diagram
\[\xymatrix{
  & & \cone (Rf) \ar[r] & RE'[1] \\
  & \widetilde{\cone f} \ar@{-->}[ur] \ar[d] & \cone f \ar[u] \ar[l] \\
  QE'' \ar@{-->}[ur] \ar[r] & E''
}\]
commutative, where the $\cone f \to \widetilde{\cone f}$ is a
cofibration and $\widetilde{\cone f} \to E''$ is an acyclic fibration. 
The extension class associated to the short exact sequence above is
given by the composition $\tilde \alpha$ of the formal inverse of $QE''
\to E''$, the composition $QE'' \to RE'[1]$ of the dashed arrows with the
morphism $\cone (Rf) \to RE'[1]$ and the formal inverse of $E'[1] \to
RE'[1]$. By abuse of notion we often call the morphism $\alpha: QE''
\to RE'[1]$ the extension class associated to the short exact sequence
above. 

\subsection{Definition of the Atiyah class}

Let $\cat C$ be a symmetric monoidal model category of cochain complexes. 

Let $\sgm O$ be an operad in $\cat C$ and let $A$ be an $\sgm
O$-algebra. Let $d: \sgm O \to M$ be a derivation where $M$ is
supposed to be a cofibrant $A$-module. We make this assumption (which
can be seen as a kind of ``smoothness'' assumption), to
simplify things in what follows.

Let $E$ be any cofibrant $A$-module. Recall that there is a short
exact sequence
\[\xymatrix{
  0 \ar[r] &
  P_A(M, E) \ar[r] &
  J_d E \ar[r] &
  E \ar[r] &
  0
}\]
of $A$-modules, the $d$-jet module sequence of $E$. 
\begin{definition}
  The extension class $\alpha_E: E \to P_A(M, E)[1]$ of the $d$-jet
  module sequence of $E$ is the \emph{(operadic) $d$-Atiyah class of
    $E$}. 
\end{definition}

\begin{remark}
  By the considerations of the previous subsection, the extension
  class is represented by a morphism $E \to RP_A(M, E)[1]$ in $\cat C$
  composed with the formal inverse of the natural morphism $P_A(M,
  E)[1] \to RP_A(M, E)[1]$. 
\end{remark}

\begin{example}
  Let $A$ be a unital, commutative algebra in $\cat C$, which we also
  view as $\Com$-algebra. Then the definition of the operadic Atiyah
  class above corresponds to the usual definition of the Atiyah class
  of a cofibrant module over a unital commutative algebra. 
\end{example}
This example is considered in~\cite{Kapranov99} where $\cat C$ is the
category of $\mathcal O_X$-modules over a smooth complex manifold $X$.

\subsection{A stabilised category}

Let $\cat C$ be a symmetric monoidal model category of cochain
complexes. Let $R$ be the symmetric monoidal fibrant replacement
functor of $\cat C^*$. 

Let $X$ and $Y$ be two objects in $\cat C^*$. Thanks to the natural
transformations $R^n Y \to R^{n + 1} Y$, there is a natural sequence
\[\hom(X, Y) \longrightarrow \hom(X, RY) \longrightarrow \hom(X, R^2
Y) \longrightarrow \dots\] of sets. In
particular, we can form the colimit $\colim_{n \in \mathbf N_0}
\hom(X, R^n Y)$.

As the morphisms $R^n Y \to R^{n + 1} Y$ are weak equivalences,
i.e.~induce isomorphisms in the homotopy category, there is a
well-defined natural map $$\colim_{n \in \mathbf N_0} \hom(X, R^n Y)
\to \Ho \cat C(X, Y).$$

\begin{definition}
  An element in the colimit $\colim_{n \in \mathbf N_0}
  \hom(X, R^n Y)$ is a \emph{stable morphism from $X$ to $Y$}. Stable
  morphisms are denoted by arrows of the form $X \stableto Y$. 
\end{definition}

\begin{remark}
  There is a canonical way to compose two stable morphisms $X
  \stableto Y$ and $Y \stableto Z$ to a stable morphism $X \stableto
  Z$. Thus, one may form the category of stable morphisms of $\cat
  C^*$ that has the same objects as $\cat C^*$
\end{remark}

Let $\sgm O$ be an operad in $\cat C^*$, and let $A$ be a $\sgm
O$-algebra. The induced fibrant replacement functor on the model
category of $A$-modules is also denoted by $R$. 
\begin{remark}
  Using this functor, one can also define stable morphisms of
  $A$-modules and the category of stable morphisms of $A$-modules.
  There is a faithful forgetful functor from the category of stable
  morphisms of $A$-modules to the category of stable morphisms of
  $\cat C^*$.
\end{remark}

\subsection{Free connections}

Let $\cat C$ be a symmetric monoidal model category of cochain
complexes. 

Let $\sgm O$ be an operad in $\cat C$ and let $A$ be an $\sgm
O$-algebra. Let $d: \sgm O \to M$ be a derivation where $M$ is
supposed to be a cofibrant $A$-module. 

Let $E$ be any cofibrant $A$-module. A stable free morphism $\nabla: E
\freestableto P_A(M, E)$ of cochain complexes over $\cat C$ is called
\emph{freely $d$-derivative} if the canonical diagrams
\[\xymatrix@C+96pt{\sgm O(n + 1) \otimes A^n \otimes E
  \ar@^{->}[r]^{\substack{\id \otimes (\sum_{p + 1 + q = n} \id^{\otimes p} \otimes d
    \otimes \id^{\otimes q}) \otimes \id \\ + \id \otimes \id^{\otimes n}
    \otimes \nabla}}_>0 \ar[d] &
  {\begin{aligned}
      \sgm O(n + 1) \otimes A^{\otimes (n - 1)} \otimes M \otimes E
      \\
      \oplus \sgm O(n + 1) \otimes A^{\otimes n} \otimes P_A(M, E)
    \end{aligned}}
  \ar@^{->}[d] \\
  E \ar@^{->}[r]_{\nabla}_>0 & P_A(M, E)
}\]
commute. 
\begin{definition}
  A \emph{free $d$-connection on $E$} is a freely $d$-derivative
  stable morphism $\nabla: E \freestableto P_A(M, E)$.
\end{definition}

\begin{remark}
  Note that a free $d$-connection is in spite of its name not a
  special case of a $d$-connection.
\end{remark}

However, every $d$-connection induces naturally a free $d$-connection
as follows:
\begin{example}
  Any $d$-connection $\nabla_0: E \to P_A(M, E)$ on $E$ composed with
  the natural stable morphism $P_A(M, E) \stableto P_A(M, A)$ gives a
  free $d$-connection on $E$. 
\end{example}

Let $f: P_A(M, E) \to J_d E$ be the canonical morphism. Consider the
diagram
\[\xymatrix{
  & \cone Rf \\
  E \ar@{-->}[ur]^{\tilde \alpha} & \cone f, \ar[l] \ar[u] }\] which
is a small part of the big diagram in the previous subsection about
general extension classes. The dashed arrow $\tilde \alpha$ making the
diagram commutative exists as $E$ is by assumption already cofibrant.

Recall that up to the differential (!)
$$\cone Rf = RP_A(E, A)[1] \oplus RE \oplus RP_A(E, A).$$
Thus the
morphism $\tilde \alpha$ given above, can be decomposed as $\tilde
\alpha = (\beta, i, \nabla): E \freeto RP_A(E, A)[1] \oplus RE \oplus
RP_A(E, A)$. In particular, $\beta$ is the extension class associated
to the $d$-jet module sequence of $E$, and $i: E \to RE$ is the
natural acyclic cofibration. As $\tilde \alpha$ is a morphism of
$A$-modules, it follows that the component $\nabla$ defines a free
$d$-connection on $E$.

Thus we have proven the following proposition:
\begin{proposition}
  There exists a free $d$-connection on $E$. 
\end{proposition}
\qed

This proposition has a kind of inverse. 
\begin{proposition}
  Let $\nabla: E \freestableto P_A(M, E)$ be a free $d$-connection.
  Then $$- [\diff, \nabla] := \nabla \circ \diff - \diff \circ \nabla:
  E \stableto P_A(E)[1]$$
  is a standard representation of the
  $d$-Atiyah class of $E$.
\end{proposition}

\begin{proof}
  Assume that the stable morphism $\nabla$ is induced by a morphism
  $\nabla: E \freeto R^N P_A(M, E)$ for some $N \gg 0$.  Set $\beta :=
  - [\diff, \nabla]$. The morphism $\tilde \alpha = (\beta, i,
  \nabla): E \to \cone R^N f$, where $i: E \to R^N E$ is the natural
  acyclic cofibration, is a morphism of $A$-modules and makes the
  diagram
  \[\xymatrix{
    & \cone R^N f \\
    E \ar[ur]^{\tilde \alpha} & \cone f \ar[l] \ar[u],
  }\]
  commutative. By definition of the extension class of the $d$-jet module
  sequence, the proposition follows. 
\end{proof}

\subsection{Derivations of morphisms}

Let $\cat C$ be a symmetric monoidal model category of cochain
complexes. 

Let $\sgm O$ be an operad in $\cat C$ and let $A$ be an $\sgm
O$-algebra. Let $d: A \to M$ be a derivation where $M$ is
supposed to be a cofibrant $A$-module. 

Let $E_1, \dots, E_m$ be cofibrant $A$-modules. Let $\nabla_i: E_i \freestableto
P_A(M, E_i)$ be free $d$-derivations. There is a free stable morphism
$$\nabla: P_A(E_1, \dots, E_m) \freestableto P_A(M, E_1, \dots,
E_m)$$
of cochain complexes over $\cat C$ which is induced by the
canonical morphisms
\[
\xymatrix{
  \sgm O(n + m) \otimes A^{\otimes n} \otimes E_1 \otimes \dots \otimes
  E_m \ar@^{->}[d]^{\sum_{j = 1}^m \nabla_i}_>0 \\
  \oplus_{j = 1}^m P_A(E_1, \dots, E_{j - 1}, P_A(M, E_j), E_{j + 1},
  \dots, E_m) \ar@^{->}[d]_>0 \\
  P_A(M, E_1, \dots, E_m)
}\]

\begin{definition}
  The free stable morphism $$\nabla: P_A(E_1, \dots, E_m) \freestableto
  P_A(M, E_1, \dots, E_m)$$ of cochain complexes is denoted by
  $(\nabla_1, \dots, \nabla_m)$.
\end{definition}

\begin{remark}
  It is easy to see that $\nabla$ as defined above, behaves similar to
  a free $d$-connection, i.e.~it makes the natural diagram 
  \[
  \xymatrix@C+120pt{
    {\begin{aligned}
        \sgm O(n + 1) \otimes A^n \\ \otimes P_A(E_1, \dots, E_m)
      \end{aligned}}
    \ar@^{->}[r]^{\substack{\id \otimes (\sum_{p + 1 + q = n}
        \id^{\otimes p} \otimes d
        \otimes \id^{\otimes q}) \otimes \id \\
        + \id \otimes \id^{\otimes n}} \otimes \nabla}_>0 \ar[d] &
    {\begin{aligned} \sgm O(n + 1) \otimes A^{\otimes (n - 1)} \\ \otimes
        M \otimes
        P_A(E_1, \dots, E_m) \\
        \oplus \sgm O(n + 1) \otimes A^{\otimes n} \\ \otimes P_A(M, E_1,
        \dots, E_m)
      \end{aligned}}
    \ar@^{->}[d]
    \\
    P_A(E_1, \dots, E_m) \ar@^{->}[r]_{\nabla}_>0 & P_A(M, E_1, \dots, E_m)
  }\]
  commutative. 
\end{remark}

Let $E$ be another $A$-module and $\nabla: E \freestableto P_A(M, E)$ a free
$d$-connection. 

Let $f: E \stableto P(E_1, \dots, E_m)$ be a stable morphism of
$A$-modules. Note that $f$ induces naturally a stable morphism $P_A(\id, f):
P_A(M, E) \stableto P_A(M, E_1, \dots, E_n)$ of $A$-modules. 

Define the free stable morphism
$$\nabla f := [\nabla, f] := (\nabla_1, \dots, \nabla_m) \circ f -
(\id, f) \circ \nabla: E \freestableto P_A(M, E_1, \dots, E_m).$$ 

\begin{definition}
  The morphism $\nabla f$ is the \emph{derivative of $f$ (with respect
    to $\nabla$ and $\nabla_1, \dots, \nabla_m$)}. 
\end{definition}

\subsection{The free $A$-algebra}

Now we are ready to come to the central subsections. The main object
of our studies will be the free algebra over a module and derivations
on it.

Let $\cat C$ be a symmetric monoidal model category of cochain
complexes.  Let $\sgm O$ be an operad in $\cat C$ and let $A$ be an
$\cat O$-algebra.

Let $d: A \to M$ be a derivation where $M$ is supposed to
be a cofibrant $A$-module. 

Let $\nabla: M \freestableto P_A(M, M)$ be a free
$d$-derivation. As detailed in the previous subsection, it induces
free stable morphisms $\nabla: P_A^n(M) \freestableto P_A^{n + 1}(M)$, $n \ge
0$. 

It is easy to see that there exist unique free stable morphisms $\nabla:
S_A^n(M) \freestableto S_A^{n + 1}(M)$ such that the canonical diagrams
\[\xymatrix{
  P_A^n(M) \ar@^{->}[r]^\nabla_>0 \ar[d] & P_A^{n + 1}(M) \ar[d] \\
  S_A^n(M) \ar@^{->}[r]_\nabla_>0 & S_A^{n + 1}(M)
}\]
commute. Here, the vertical arrows denote the canonical projections
onto the coinvariants. 

A free stable morphism $Q: S_A^*(M) \freestableto S_A^*(M)$ is
\emph{freely derivative over $d$} if $Q$ is a free stable derivation
and the diagram
\[\xymatrix{
  A \ar[r]^d \ar[d] & M \ar[d] \\
  S_A^*(M) \ar@^{->}[r]_Q_>0 & S_A^*(M)
}\]
commutes. 

\begin{definition}
  A free stable morphism $Q: S_A^*(M) \freestableto S_A^*(M)$ that is
  freely derivative over $d$ is a
  \emph{free derivation over $d$}.
\end{definition}

An example of such a free derivation over $d$ is given as follows.
First note that the free stable morphisms $\nabla: S^n_A(M)
\freestableto S^{n + 1}_A(M)$ as constructed above altogether induce a
free stable morphism $Q_\nabla: S^*_A(M) \freestableto S^*_A(M)$. It
easily follows that:
\begin{proposition}
  The morphism $Q_\nabla: S_A^*(M) \freestableto S_A^*(M)$ is a free
  derivation over $d$. 
\end{proposition}
\qed

Let $E$ be another cofibrant $A$-module. Let $\nabla_E: E \freestableto
P_A(M, E)$ be a free $d$-connection. 

By analoguous considerations as above, there exist unique free stable
morphisms $\nabla_E: S_A^n(M, E) \freestableto S_A^{n + 1}(M, E)$ such
that the canonical diagrams
\[\xymatrix{
  P_A^n(M, E) \ar@^{->}[r]^{\nabla_E}_>0 \ar[d] & P_A^{n + 1}(M, E) \ar[d] \\
  S_A^n(M, E) \ar@^{->}[r]_{\nabla_E}_>0 & S_A^{n + 1}(M, E)
}\]
commute. Again, the vertical arrows denote the canonical projections
onto the coinvariants. 

Let $Q: S_A^*(M) \freestableto S_A^*(M)$ be a free derivation over
$d$.  A free stable morphism $D: S_A^*(M, E) \freestableto S_A^*(M,
E)$ is \emph{freely $Q$-derivative over $\nabla_E$} if $D$ is a free
$Q$-connection and the canonical diagram
\[\xymatrix{
  E \ar@^{->}[r]^{\nabla_E}_>0 \ar[d] & P_A(M, E) \ar@^{->}[d]_>0 \\
  S^*_A(M, E) \ar@^{->}[r]_D_>0 & S^*_A(M, E)
}\]
commutes. 
\begin{definition}
  A free stable morphism $D: S_A^*(M, E) \freestableto S_A^*(M, E)$ that is
  freely $Q$-derivative over $\nabla_E$ is a
  \emph{free $Q$-derivation over $\nabla_E$}.
\end{definition}

Again we can easily construct an example: By $D_\nabla$ we denote the
free stable morphism $S_A^*(M, E) \freestableto S_A^*(M, E)$ induced by the
morphisms $\nabla_E: S^n(M, E) \freestableto S^{n + 1}(M, E)$.
\begin{proposition}
  The free stable morphism $D_\nabla: S_A^*(M, E) \freestableto
  S_A^*(M, E)$ is a free $Q_\nabla$-connection over $\nabla_E: E
  \freestableto P_A(M, E)$.
\end{proposition}

\subsection{The total curvature form}

Let $\cat C$ be a symmetric monoidal model category of cochain
complexes. We assume that $\cat C$ is $\set Q$-linear, i.e.~that the
hom-sets are $\set Q$-vector spaces.

Let $\sgm O$ be an operad in $\cat C$ and let $A$ be an $\sgm
O$-algebra. Let $d: \sgm O \to M$ be a derivation where $M$ is
supposed to be a cofibrant $A$-module. 

Let $\nabla: M \freestableto P_A(M, M)$ be a free stable
$d$-connection. Set
\[
R_\nabla := \exp(\ad Q_\nabla) \diff = \sum_{n = 0}^\infty \frac 1 {n!} 
(\ad Q_\nabla)^n(\diff): S_A^*(M) \stableto S_A^*(M)[1],
\]
where, as always, $\diff: S_A^*(M) \to S_A^*(M)[1]$ denotes the differential of
$S_A^*(M)$. We note that $R_\nabla$ is a stable morphism and a stable
derivation of $S_A^*(M)$ of degree one over $d$. 
\begin{definition}
  The stable morphism $R_\nabla: S_A^*(M) \stableto S_A^*(M)[1]$ is
  the \emph{total curvature form of $M$ associated to $\nabla$}.
\end{definition}

In fact, $R_\nabla$ is a differential:
\begin{proposition}
  One has $R_\nabla \circ R_\nabla = \frac 1 2 [R_\nabla, R_\nabla] =
  0$. 
\end{proposition}

\begin{proof}
  This follows from \[[R_\nabla, R_\nabla] = [\exp(\ad Q_\nabla)
  \diff, \exp(\ad Q_\nabla) \diff] = \exp(\ad Q_\nabla) [\diff, \diff]
  = 0.\]
\end{proof}

Denote by $R^{(n)}_\nabla: M \stableto S_A^n(M)[1]$, $n \ge 0$ be those
stable morphisms such that the diagrams
\[\xymatrix{
  M \ar[d] \ar@^{->}[r]^{R_\nabla^{(n)}} & S_A^n(M) \ar@^{->}[d] \\
  S_A^*(M) \ar@^{->}[r]_{R_\nabla} & S_A^*(M)
}\]
commute. 

\begin{example}
  We have $R_\nabla^{(0)} = 0$ and $R_\nabla^{(1)} = \diff$, the
  differential of the complex $M$. Furthermore, we have
  \[
  R_\nabla^{(2)} = -[\diff, \nabla]: M \stableto S_A^2(M),\]
  i.e.~$R_\nabla^{(2)}$ represents the symmetrisation of the Atiyah
  class of $M$ in the homotopy category of $A$-modules.
\end{example}

Let $E$ be a cofibrant $A$-module and $\nabla_E: E \freestableto
P_A(M, E)$ be a free stable $d$-connection. Set
\[
T_\nabla := \exp(\ad D_\nabla) \diff = \sum_{n = 0}^\infty \frac 1 {n!} 
(\ad D_\nabla)^n(\diff): S_A^*(M, E) \stableto S_A^*(M, E)[1]
\]
where $\partial: S_A^*(M, E) \to S_A^*(M, E)[1]$ denotes the
differential of $S_A^*(M)$. We note that $T_\nabla$ is a stable
derivation of $S_A^*(M, E)$ of degree one over $d$.
\begin{definition}
  The morphism $T_\nabla: S_A^*(M, E) \to S_A^*(M, E)[1]$ is the \emph{total
    curvature form of $E$ associated to $\nabla$ and $\nabla_E$}. 
\end{definition}

In fact, $T_\nabla$ is a differential. The proof goes as in the case
of $R_\nabla$.
\begin{proposition}
  One has $T_\nabla \circ T_\nabla = \frac 1 2 [T_\nabla, T_\nabla] =
  0$. 
\end{proposition}
\qed

Denote by $T^{(n)}_\nabla: M \stableto S_A^n(M, E)[1]$, $n \ge
0$ those stable morphisms such that the diagrams
\[\xymatrix{
  E \ar[d] \ar@^{->}[r]^{T_\nabla^{(n)}} & S_A^n(M, E) \ar@^{->}[d] \\
  S_A^*(M, E) \ar@^{->}[r]_{T_\nabla} & S_A^*(M, E)
}\]
commute. 

\begin{example}
  We have $T_\nabla^{(0)} = \diff$, the differential of the complex $M$. 
  Furthermore, we have
  \[
  T_\nabla^{(1)} = -[\diff, \nabla]: M \stableto P_A(M, E),\]
  i.e.~$T_\nabla^{(1)}$ represents the Atiyah class of $E$ in the
  homotopy category of $A$-modules.
\end{example}

We owe the idea to define the total curvature classes $R_\nabla$ and
$T_\nabla$ the article~\cite{Kapranov99} by Kapranov, where this is
done in the classical geometrical context.

\subsection{Bianchi identity}

Let $\cat C$ be a $\set Q$-linear symmetric monoidal model category of
cochain complexes.

Let $\sgm O$ be an operad in $\cat C$ and let $A$ be an $\sgm
O$-algebra. Let $d: \sgm O \to M$ be a derivation where $M$ is
supposed to be a cofibrant $A$-module. 

Let $\nabla: M \freestableto P_A(M, M)$ be a free stable
$d$-connection. Let $\alpha: M \to S^2_A(M)[1]$ be the (symmetrised)
Atiyah class of $M$ (i.e.~it is a morphism in the homotopy
category). It induces a morphism $\hat \alpha: S^2_A(M) \to
S^3_A(M)[1]$ such that the canonical diagram (in the homotopy category)
\[
\xymatrix@C+60pt{ P_A(M, M) \ar[r]^{P_A(\alpha, \id) + P_A(\id,
    \alpha)} \ar[d] &
  P_A(M, M, M)[1] \ar[d] \\
  S_A^2(M) \ar[r]_{\hat \alpha} & S_A^3(M)[1] }\] commutes.

The following proposition is known under the name ``homological Bianchi
identity'' in the classical, i.e.~non-operadic,
case. See~\cite{Kapranov99}.
\begin{proposition}
  It is
  \[
  \hat \alpha \circ \alpha: M \to S^3_A(M)[2]
  \]
  the zero morphism in the homotopy category. 

  This is called the Bianchi identity of the Atiyah class. 
\end{proposition}

\begin{proof}
  A representative of $\alpha$ is the free stable morphism
  $R_\nabla^{(2)}: M \freestableto S^2_A(M)$. Let $\hat
  R_\nabla^{(2)}: S^2(M) \freestableto S^3_A(M)$ be the free stable
  morphism such that the canonical diagram
  \[
  \xymatrix@C+72pt{ P_A(M, M) \ar@^{->}[r]^{P_A(R_\nabla^{(2)}, \id) +
      P_A(\id, R_\nabla^{(2)})}_>0 \ar[d] &
    {\begin{aligned}
        P_A(S^2_A(M), M) \\
        \oplus P_A(M, S^2_A(M))
      \end{aligned}}
    \ar@^{->}[d]_>0 \\
    S_A^2(M) \ar@^{->}[r]_{\hat R_\nabla^{(2)}}_>0 & S_A^3(M) }\]
  commutes, i.e.~$\hat R_\nabla^{(2)}$ is a representative of $\hat
  \alpha$.
  
  From $[R_\nabla, R_\nabla] = 0$ we deduce $\hat R_\nabla^{(2)} \circ
  R_\nabla^{(2)} = [\partial, R_\nabla^{(3)}]$ by looking at the terms
  of low degree. Thus the left hand side, which gives $\hat \alpha
  \circ \alpha$ in the homotopy category, is null-homotopic, which
  means that it is the zero morphism in the homotopy category.
\end{proof}

We can do the same for the Atiyah class of a cofibrant $A$-module $E$.
Let $\nabla_E: M \freestableto P(M, E)$ be a free stable $d$-connection. Let $\alpha_E: M
\to P_A(M, E)$ be the Atiyah class of $E$. Together with (the
non-symmetrised version of) $\alpha$ it is induces a morphism $\hat
\alpha_E: S^1(M, E) \to S^2_A(M, E)[1]$ such that the canonical diagram
\[
\xymatrix@C+60pt{ P_A(M, E) \ar[r]^{P_A(\alpha, \id) + P_A(\id,
    \alpha_E)} \ar[d] &
  P_A(M, M, E)[1] \ar[d] \\
  S_A^1(M, E) \ar[r]_{\hat \alpha_E} & S_A^2(M, E)[1]
}\]
commutes. 

Again we have a cohomological Bianchi identity:
\begin{proposition}
  It is
  \[
  \hat \alpha_E \circ \alpha_E: M \to S^2_A(M, E)[2]
  \]
  the zero morphism in the homotopy category. 
\end{proposition}
\qed

\subsection{An example}

Let $\cat C$ be a $\set Q$-linear symmetric monoidal model category of
cochain complexes.

Let $\sgm O$ be an operad in $\cat C$ and let $V$ be a cofibrant
object in $\cat C$. Consider the free $\sgm O$-algebra $A := F_{\sgm
  O}(V)$. 

Recall that any free morphism $g: V \freeto A[1]$ induces a free
derivation $\hat g: A \to A[1]$.

\begin{definition}
  The morphism $g$ is a \emph{solution of the Maurer--Cartan equation
    for $A$} if and only if
  \[
  [\diff + \hat g, \diff + \hat g] = 0,
  \]
  where $\diff$ is the differential on $A$.
\end{definition}

In other words, $g$ is a solution if and only if $\diff + \hat g$
defines a new differential on $A$. In that case, we denote by $A(g)$
the $\sgm O$-algebra we get when we substitute its differential $\diff$ by
$\diff + \hat g$. We consider $A(g)$ as a \emph{deformation of
  $A$}. It is not a free algebra anymore.

An example of such a solution is given by a (strong homotopy) Lie
coalgebra structure on $V[1]$, see, e.g.,~\cite{Kontsevich99}.

\begin{definition}
  Two solutions $g_0$ and $g_1$ of the Maurer--Cartan equation for $A$
  are \emph{gauge equivalent} if and only if they appear simultaneouly
  in a family $g(t)$ of solutions with
  \[
  \hat g'(t) = [\hat \xi(t), \diff + \hat g(t)]
  \]
  for a family $\xi(t)$ of morphisms $V \to A$ (that induces a
  family $\hat \xi(t)$ of derivations $A \to A$). 
\end{definition}

The free $A$-module $M := F_{A}(V)$ is a cofibrant $\sgm
O$-module by definition of the model structure on the category of
$A$-modules as $V$ is cofibrant.

Denote by $d: A \to M$ the derivation of $A$ into $M$ that is induced
by the identity map $V \to V$. If $g$ is a solution of the
Maurer--Cartan equation for $A$, we denote by $M(g)$ the object $M$
with the unique differential $\diff(g)$ such that it becomes an
$A(g)$-module and $d$ a morphism from $A(g)$ to $M(g)$.

Recall further that we have a canonical $d$-connection $\nabla: M \to
P_{A}(M, M)$. This is a free $d$-connection of $A(g)$-modules when
viewed as a morphism $\nabla: M(g) \freeto P_{A(g)}(M(g), M(g))$.

Thus we can consider its total curvature form $R(g)_\nabla = \exp(\ad
\nabla) (\diff(g))$. Set $\hat R(g) := R_\nabla (g) - \diff_0$,
where $\diff_0$ on the right hand side is the differential on
$S_A(M)$. Let $R(g): M(g) \freeto S_A(M(g))$ be the restriction of $\hat
R(g)$ to $M(g)$ considered as a free morphism of $A(g)$-modules. As
$[R(g)_\nabla, R(g)_\nabla] = 0$, we have the following proposition:
\begin{proposition}
  The morphism $R(g): M(g) \to S_A(M(g))$ viewed a free morphism
  $R(g): M \to S_A(M)$ is a solution of the Maurer--Cartan equation
  for $S_A(M)$, i.e.~$[\diff_0 + \hat R(g), \diff_0 + \hat R(g)] = 0$.
\end{proposition}

In fact, we have the following:
\begin{proposition}
  Let $g_1$ and $g_2$ be gauge equivalent solutions of the
  Maurer--Cartan equation for $A$. Then $R(g_1)$ and $R(g_2)$ are gauge
  equivalent. 
\end{proposition}

\begin{proof}
  Let $g(t)$ be a family of solutions connecting $g_1$ and $g_2$ with
  $\hat g'(t) = [\hat \xi(t), \diff + \hat g(t)]$. Then
  $\hat R(g(t)) = [\exp(\ad \nabla) \hat \xi(t), \diff_0 + \hat R(g(t))]$. 
\end{proof}

\begin{remark}
  Thus, what we have defined is a canonical map from the set of
  solutions of the Maurer--Cartan equation for $A$ up to gauge
  equivalence to set of solutions of the Maurer--Cartan equation for
  $S^*_A(M)$ up to gauge equivalence. 
\end{remark}

\appendix

\section{Proofs}

\subsection{Smallness in Grothendieck categories}

In this appendix we give a simple proof for the following theorem:
\begin{theorem}
  Any object in a Grothendieck category is small. 
\end{theorem}

\begin{proof}
  Let $A$ be an object in a Grothendieck category $\cat C$. Let $G$ be
  a generator of $\cat C$. Let $\mathfrak A$ be the set of subobjects
  of $A$ and set $\hat G := \bigoplus_{A' \in \mathfrak A} G$. Let
  $\kappa$ be the sum of the cardinality of the set $\mathfrak A$ of
  subobjects of $A$ and the cardinality of the set of subobjects of
  $\hat G$. We show that $A$ is $\kappa$-small. Let $\lambda$ be an
  ordinal with $\cofin \lambda > \kappa$ and let $X: \lambda \to
  \mathcal C$ be a colimit-preserving functor. We have to show that
  the natural map $\colim_{\mu < \lambda} \hom(A, X_\mu) \to \hom(A,
  \colim_{\mu < \lambda} X_\mu)$ is an isomorphism.
  
  First, we show the surjectivity. Let $f: A \to X_\lambda :=
  \colim_{\mu < \lambda} X_\mu$ be any morphism in $\mathcal C$. We
  have to show that $A$ factors through some $X_\beta$ for $\beta <
  \lambda$. By dividing out the kernel if necessary, we may assume
  that $f$ is a monomorphism. For each subobject $A'$ of $A$ there is
  a $\beta < \lambda$ such that the map from the fibre product
  $X_\beta \times_{X_\lambda} A \to X_\lambda \to X_\lambda/f(A')$ is
  not the zero map if $A'$ is a proper subobject of $A$. As the
  cofinality of $\lambda$ is greater than the cardinality of the set
  of subobjects of $A$, there is a $\beta < \lambda$ such that
  $X_\beta \times_{X_\lambda} A \to X_\lambda \to X_\lambda/f(A')$ is
  not the zero map for each proper subobject $A'$ of $A$.  For each
  subobject $A'$ there is a morphism $G \to A$ and a morphism $G \to
  X_\gamma$ such that $G \to A \to X_\lambda$ equals $G \to X_\gamma
  \to X_\lambda$ and such that $G \to A \to X_\lambda \to
  X_\lambda/A'$ is non-zero if $A'$ is a proper subobject of $A$.
  These morphisms induce by the universal property of the direct sum
  $\hat G$ a morphism $\hat G \to A$ and a morphism $\hat G \to
  X_\beta$ such that $\hat G \to A \to X_\lambda$ equals $\hat G \to
  X_\beta \to X_\lambda$.  Let $A'$ be the image of the morphism $\hat
  G \to A$. Suppose that $A'$ is a proper subobject of $A$. Then $A
  \to X_\lambda \to X_\lambda/f(A')$ is not the zero map. In
  particular, $\hat G \to A \to X_\lambda \to X_\lambda/f(A')$ is not
  the zero map which contradicts the fact that the image of $\hat G
  \to A$ is $A'$.  Thus $A'$ cannot be a proper subobject of $A$, thus
  $A' = A$, i.e.~$\hat G$ is an epimorphism. Let us denote the kernel
  of $\hat G \to A$ be $K$. For each ordinal $\gamma$ with $\beta \leq
  \gamma < \lambda$ denote the kernel of $\hat G \to X_\beta \to
  X_\gamma$ by $K_\gamma$. As $A \to X_\lambda$ is a monomorphism, the
  kernel $K_\gamma$ is a subobject of $K$ for all ordinals $\gamma$
  with $\beta \leq \gamma < \lambda$. Consider the filtered colimit
  over $\lambda$ of the left exact sequence $0 \to K_\gamma \to \hat G
  \to X_\gamma$. As $\mathcal C$ is a Grothendieck category, this
  yields the exact sequence $0 \to \colim_{\gamma < \lambda} K_\gamma
  \to \hat G \to X_\lambda$, which shows that $\colim_{\gamma <
    \lambda} K_\gamma = K$. As the cofinality of $\lambda$ is greater
  than the cardinality of the set of subobjects of $\hat G$ the
  sequence $(K_\gamma)_{\gamma < \lambda}$ stabilises, i.e.~there is a
  $\gamma < \lambda$ with $K_\gamma = \colim_{\gamma < \lambda}
  K_\gamma = K$.  Thus the morphism $\hat G \to X_\gamma$ induces a
  morphism $\hat G/K \to X_\gamma$ and this morphism composed with the
  inverse of the isomorphism $\hat G/K \to A$ yields a morphism $A \to
  X_\gamma$ such that $A \to X_\gamma \to X_\lambda$ equals $f$.
  
  We now show the injectivity. Let $A \to X_{\beta_1}$ and $A \to
  X_{\beta_2}$ with $\beta_1, \beta_2 < \lambda$ be two morphisms such
  that $A \to X_{\beta_1} \to X_\lambda$ and $A \to X_{\beta_2} \to
  X_\lambda$ are equal. We have to show that there is an ordinal
  $\gamma$ with $\beta_1, \beta_2 \leq \gamma < \lambda$ such that $A
  \to X_{\beta_1} \to X_\gamma$ and $A \to X_{\beta_2} \to X_\gamma$
  are equal. For this, let $K_\gamma$ be the equaliser of $A \to
  X_{\beta_1} \to X_\gamma$ and $A \to X_{\beta_2} \to X_\gamma$ for
  any ordinal $\gamma$ with $\beta_1, \beta_2 \leq \gamma < \lambda$.
  Consider the filtered colimit over the exact sequence $0 \to
  K_\gamma \to A \rightrightarrows X_\gamma$. As $\mathcal C$ is a
  Grothendieck category, this yields the exact sequence $0 \to
  \colim_{\gamma < \lambda} K_\gamma \to A \rightrightarrows
  X_\lambda$, which shows that $\colim_{\gamma < \lambda} K_\gamma =
  0$. As the cofinality of $\lambda$ is greater than the set of
  subobjects of $A$, the sequence $(K_\gamma)_{\gamma < \lambda}$
  stabilises. In particular, there is a $\gamma < \lambda$ with
  $K_\gamma = \colim_{\delta < \lambda} K_\delta = 0$. Thus $A \to
  X_{\beta_1} \to X_\gamma$ and $A \to X_{\beta_2} \to X_\gamma$ are
  equal.
\end{proof}

\bibliographystyle{amsplain}
\bibliography{jets_and_operads}

\end{document}